\documentclass[12pt]{amsart}

\usepackage{color}

\usepackage{amssymb}
\usepackage{graphicx}  

\newcommand{\CC}{{\mathbb C}}

\newcommand{\ZZ}{{\mathbb Z}}

\newcommand{\RR}{{\mathbb R}}

\newcommand{\defeq}{\stackrel{\rm{def}}{=}}

\newcommand{\rest}{\!\!\restriction}

\renewcommand{\Re}{\mathop{\rm Re}\nolimits}
\renewcommand{\Im}{\mathop{\rm Im}\nolimits}

\theoremstyle{plain}

\newtheorem{thm}{Theorem}
\newtheorem{prop}{Proposition}[section]

\newtheorem{lem}[prop]{Lemma}

\theoremstyle{definition}

\numberwithin{equation}{section}

\def\squarebox#1{\hbox to #1{\hfill\vbox to #1{\vfill}}} 
\newcommand{\stopthm}{\hfill\hfill\vbox{\hrule\hbox{\vrule\squarebox 
                 {.667em}\vrule}\hrule}\smallskip} 
\newcommand{\sech}{\textnormal{sech}}
\newcommand{\indentalign}{\hspace{0.3in}&\hspace{-0.3in}}
\newcommand{\la}{\langle}
\newcommand{\ra}{\rangle}

\usepackage{amsxtra}

\ifx\pdfoutput\undefined
  \DeclareGraphicsExtensions{.pstex, .eps}
\else
  \ifx\pdfoutput\relax
    \DeclareGraphicsExtensions{.pstex, .eps}
  \else
    \ifnum\pdfoutput>0
      \DeclareGraphicsExtensions{.pdf}
    \else
      \DeclareGraphicsExtensions{.pstex, .eps}
    \fi
  \fi
\fi

\title
[Slow soliton interaction with delta impurities]
{Slow soliton interaction with delta impurities}

\author[J. Holmer]
{Justin Holmer}
\email{holmer@math.berkeley.edu}
\author[M. Zworski]
{Maciej Zworski}
\email{zworski@math.berkeley.edu}
\address{Mathematics Department, University of California \\
Evans Hall, Berkeley, CA 94720, USA}

\begin{document}

\begin{abstract}
We study the Gross-Pitaevskii equation with a delta function potential, 
$ q \delta_0 $, where $ |q| $ is small
and analyze the solutions for which the initial condition is a 
soliton with initial velocity $ v_0 $.
We show that up to time $ (|q| + v_0^2 )^{-\frac12} \log(1/|q|) $ 
the bulk of the solution is a soliton evolving 
according the classical dynamics of a natural 
effective Hamiltonian, $ (\xi^2 + q \, \sech^2 ( x ) )/2 $.

\end{abstract}

\maketitle

\section{Introduction}   
\label{in}

The Gross-Pitaevskii equation (NLS) with a 
delta function potential and soliton initial data, 
\begin{equation}
\label{eq:nls}
\left\{
\begin{aligned}
&i\partial_t u + \tfrac{1}{2}\partial_x^2 u -q\delta_0(x)u +u|u|^2 = 0\\
&u(x,0) = e^{ i v_0 x } \sech ( x - a_0 ) \,,  
\end{aligned}
\right.
\end{equation}
offers a surprising wealth of dynamical phenomena. In \cite{HMZ1},
(and numerically in \cite{HMZ2}), the authors and 
J. Marzuola studied the high velocity, $ v_0 \gg 1 $, case and showed
that the scattering matrix of the delta potential controls the
dynamics. In this paper we describe the case of small $ q $. The
most interesting dynamics is visible for 
initial  velocities satisfying $ v_0^2 \leq |q| $. 
The low $ v_0  $  regime has been studied in the physics 
literature \cite{CM},\cite{GHW},\cite{BL}, and the behaviour in the
intermediate
range of $ q$'s and $v_0$'s, that is 
 between the fully quantum and semiclassical cases studied
in \cite{HMZ1} and in this paper respectively, is still unclear. We state
the main result here with a slightly more precise version given in 
Theorem \ref{t:2} in 
\S \ref{pr} below.

\begin{thm}
\label{t:1}
Suppose that in \eqref{eq:nls} 
we have $  |q|\ll 1 $.
Then, on a time interval $0\leq t \leq \delta (v_0^2+|q|)^{-1/2}\log(1/|q|)$, 
\begin{equation}
\label{eq:t1}
\| u ( t ,\bullet ) -  e^{i\bullet v(t) }e^{i\gamma ( t ) } \sech
(\bullet -a ( t) ) \|_{H^1 ( \RR ) } \leq C |q|^{1-3 \delta} \,,
\end{equation}
where $ a$, $ v$, and $ \gamma$ solve the following system
of equations
\begin{equation}
\label{eq:t3}
\frac{d}{dt} {  a} =   v \,, \ \ 
\frac{d}{dt} {  v} = -\frac12 q\partial_x (\sech^2) (  a) \,, 
\ \ 
\frac{d}{dt} {  \gamma} = \frac12 + \frac {v^2} 2 
- q \sech^2 (   a ) - \frac12 q\partial_x (\sech^2) (  a) \,, 
\end{equation}
with initial data $(a_0,v_0,0)$.
\end{thm}

\begin{figure}[htp]
$$\hspace{-0.1in}\includegraphics[width=6.5in]{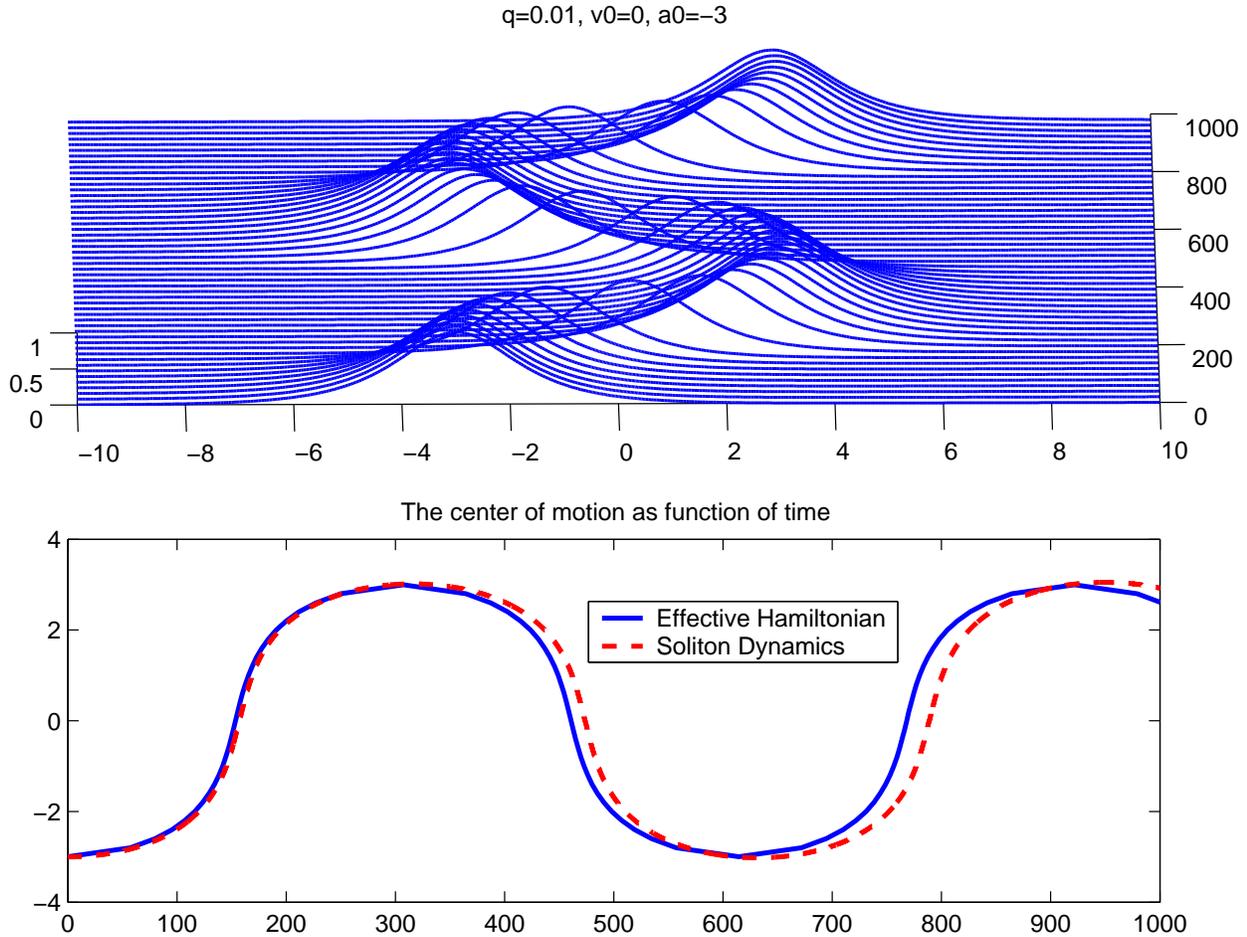}$$
\caption{The top figure shows the evolution of $ | u ( x , t ) | $ for 
 $ v_0 = 0$, $ a_0=-3$, $ q=-0.01$ for $ 0 \leq t \leq 1000$.
In the bottom figure
the dashed curve is the computed center of motion, and the
continuous curve, the plot of $   a ( t) $ given by \eqref{eq:t3}. 
More figures illustrating other cases, some with an even more dramatic
agreement can be found at 
{\tt http://math.berkeley.edu/$\sim$zworski/HZ1.pdf} }
\label{f:1}
\end{figure}

Compared to numerical results, 
the theorem gives a remarkably good description of the dynamics 
of a slow soliton interacting with a small delta function potential. 
For example consider $v_0=0$, $a_0<0$ fixed, and $|q|\to 0$, illustrated in 
Fig.\ref{f:1}. When $ q < 0 $, 
the bulk of the solution is oscillatory about the origin, with the 
center moving from $a_0<0$ to $-a_0>0$.  Since
$$\frac12 v^2 + \frac12q\eta^2(a) = \frac12q\eta^2(a_0)\,, \ \ 
\eta ( x ) \defeq \sech ( x ) \,, $$
the time to complete one cycle of oscillation is
$$\int_{a_0}^{-a_0} \frac{2\, dx}{|q|^{1/2}\sqrt{\eta^2(x) - \eta^2(a_0)}}$$
which is of size comparable to 
$ |q|^{-1/2}$.  Since the theorem provides an accurate description up to time $\sim |q|^{-1/2}\log(1/|q|)$, it covers many cycles for small enough $|q|$. When $ q > 0 $ the soliton is repulsed by the $ \delta $ potential
and slowly slides to negative infinity with the terminal velocity 
$ q^{1/2} $ -- see Fig.\ref{f:turn} below.

\renewcommand\thefootnote{\dag}%

\begin{figure}
\begin{center}
\includegraphics[width=6in]{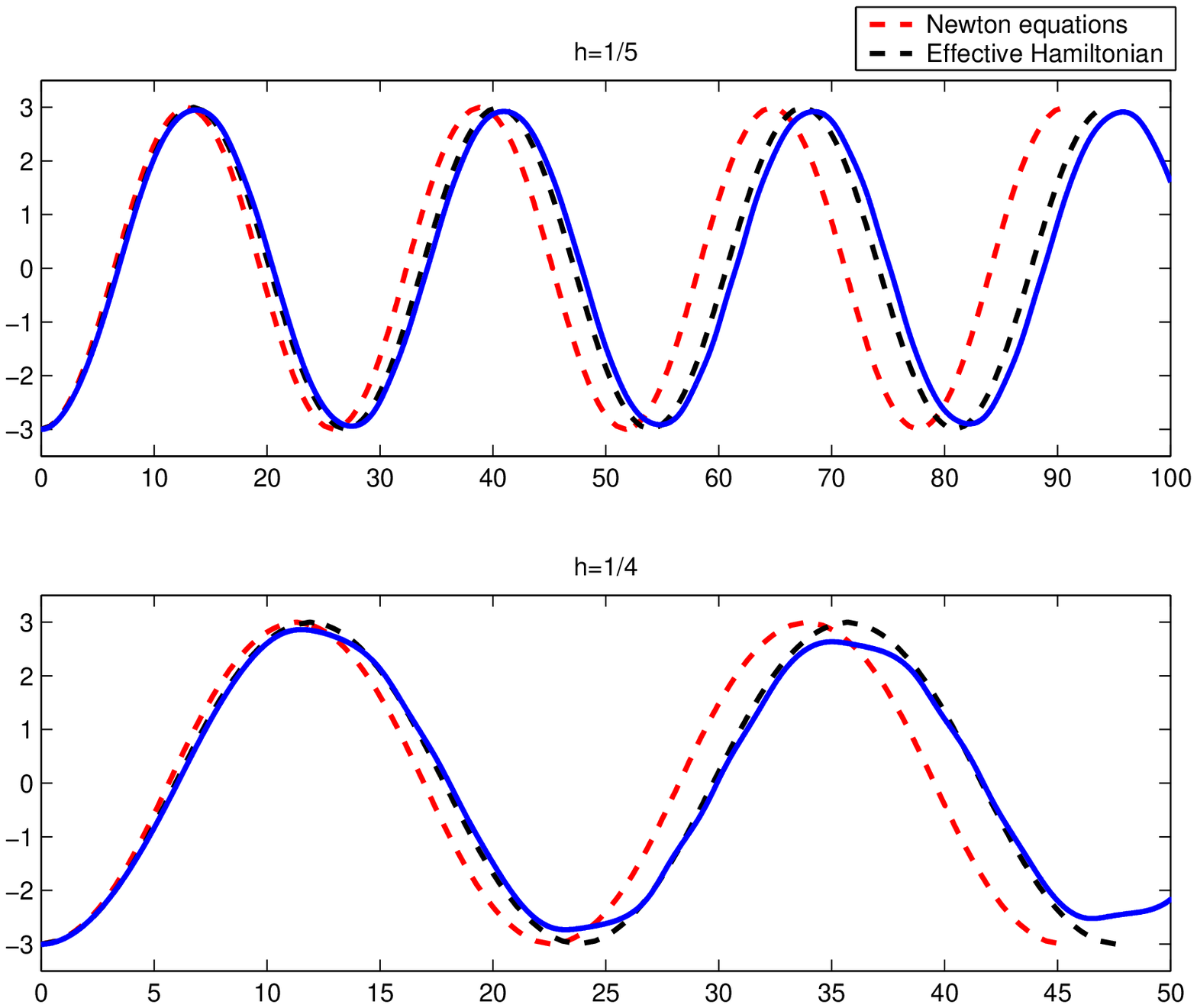}
\end{center}
\caption{Comparison of the dynamics of the center of motion of
the soliton for the Gross-Pitaevskii equation with a slowly 
varying potential, 
$$ iu_t = -\frac12 u_{xx} - |u|^2 u - \sech^2 ( h x ) u \,, \ \ 
h = 1/5 \,, \ \ h=1/4 \,, $$
and initial condition in \eqref{eq:nls} with $ v_0=0$, $ a_0=-3$.
The dashed {{red}} curve shows the solution to Newton's equations
used in \cite{BJ} and \cite{FrSi}, the {{blue}} curve shows
the center of the approximate soliton $ u $,  and the black dashed
curve is given by the equations of motion of the effective 
Hamiltonian
\[  \hspace{-1in} 
\frac 12 \left( v^2  +  \sech^2 ( h \bullet ) * \sech^2 (a )  \right) \,.\]
The improvement of the approximation given by the effective Hamiltonian
is remarkable even in the case of $ h = 1/4 $ in which we already see
radiative dissipation in the first cycle.}
\label{f:drama}
\end{figure}

The proof of our theorem follows the long tradition of the study 
of stability of solitons which started with the work of M.I. Weinstein
\cite{We}. The interaction of solitons with external potentials was 
studied in the stationary semiclassical setting by 
Floer and A. Weinstein \cite{FlWe} and Oh \cite{Oh}, and the first
dynamical result belongs to Bronski and Jerrard \cite{BJ}. 
The semiclassical 
regime is equivalent to considering slowly varying potentials,
\begin{equation}
\label{eq:slow}
\left\{
\begin{aligned}
&i\partial_t u + \tfrac{1}{2}\partial_x^2 u - W ( h x ) u +u|u|^2 = 0\,, \ \
0 < h \ll 1 \\
&u(x,0) = e^{ i v_0 x } \sech ( x - a_0 ) \,,  \ \ \| W^{ ( k )} \|_{\infty}
\leq C \,, \ \ k \leq 2 \,, 
\end{aligned}
\right.
\end{equation}
and that case 
has been studied in various settings and degrees of generality in 
\cite{FrSi}, \cite{FrSi1}, \cite{FrY} (see these papers for additional
references). The approach of these works was our starting point.
The results of \cite{FrSi} in the special case of \eqref{eq:slow}
give 
\begin{equation}
\label{eq:tslow}
\| u ( t ,\bullet ) -  e^{i\bullet v(t) }e^{i\gamma ( t ) } 
\sech(\bullet -a ( t) ) \|_{H^1 ( \RR ) } \leq C h \,,  \ \ 
0 \leq t \leq C \log(1/h)/h \,, 
\end{equation}
where 
\begin{gather}
\label{eq:t3slow}
\begin{gathered}
\frac{d}{dt} {  a} =   v + {\mathcal O} ( h^2 ) \,, \ \ 
\frac{d}{dt} {  v} =  - h W' ( h a) + {\mathcal O} ( h^2 ) \,, 
\\ 
\frac{d}{dt} {  \gamma} = \frac12 + \frac {v^2} 2 
- W ( h a ) +  {\mathcal O} ( h^2 ) \,, 
\end{gathered}
\end{gather}
with initial data $(a_0,v_0,0)$.\footnote{Strictly speaking the 
result in \cite{FrSi} describes the dynamics for $ 0 \leq t \leq c_0 / h $ 
only. That corresponds to small time dynamics of the potential $ W $.
Iterating the full strength of the result of \cite{FrSi} seems to 
give the expected extension to Ehrenfest time $ \log (1/h) /h $ \cite{Gus}.}
We note that unlike in \eqref{eq:t3} the ordinary differential system
\eqref{eq:t3slow}
is not exact -- see Fig.\ref{f:drama} and the discussion below.

At first the equations \eqref{eq:nls} and \eqref{eq:slow} appear to be very 
different: a delta function potential is very far from being slowly 
varying. The similarity of \eqref{eq:t3} and \eqref{eq:t3slow} is 
however a result
of the same underlying structure. As we recall in \S \ref{rhs} the
Gross-Pitaevski equations, \eqref{eq:nls} or \eqref{eq:slow}, are the 
equation for Hamiltonian flow of 
\begin{equation}
\label{eq:HV}
 H_V ( u ) \defeq \frac14 \int ( |\partial_x u |^2 - |u|^4 ) dx
+ \frac12 \int V  | u|^2 \,, \ \ 
V  = q \delta_0  \,, \ \ V = W ( h \bullet ) \,, 
\end{equation}
with respect to the symplectic form on $ H^1 ( \RR , \CC ) $ (considered
as a real Hilbert space):
\begin{equation}
\label{eq:ome1}
 \omega ( u , v ) = \Im \int u \bar v \,, \ \ u , v \in H^1 ( \RR , \CC) 
\,.\end{equation}
When $ V \equiv 0 $, $ \eta = \sech $ is a critical value (minimizer)
of $ H_0 $ with prescribed $ L^2 $ norm:
\begin{equation}
\label{eq:defE} d {\mathcal E}_{\eta } = 0 \,, \ \ {\mathcal E}( u ) 
\defeq H_0 ( u ) + \frac14 \| u \|_{L^2}^2 \,.
\end{equation}
The flow of $ H_0 $ is tangent to the 
manifold of solitons,
\[  M = \{ e^{i \gamma} e^{ i v ( x - a ) } \mu \, 
\sech ( \mu ( x - a ) ) \,, \ \ 
a, v, \gamma \in \RR \,, \ \ \mu \in \RR_+ \} \,, \]
which of course corresponds to the fact that the solution of \eqref{eq:nls}
with $ q = 0 $ and $ u_0 ( x , 0 ) = e^{ i \gamma + i v_0 ( x - a_0 ) } \mu
 \sech ( \mu ( x - a_0 ) ) $, is
\begin{equation}
\label{eq:frees}  u ( x , t) = 
e^{ i \gamma + i v_0  ( x - a_0) + i ( \mu^2 - v^2) t/ 2 } 
\mu \, \sech ( \mu ( x - a_0 - v_0 t ) )
\,. \end{equation}
The symplectic form \eqref{eq:ome1} restricted to $ M $ is
\begin{equation}
\label{eq:ome2} \omega\rest_M
 = \mu dv\wedge da + vd\mu\wedge da + d\gamma \wedge d\mu
\,,\end{equation}
see \S \ref{ssms}. The evolution of the parameters $ ( a, v ,\gamma, \mu ) $
in the solution $ u ( x, t ) $ follows the Hamilton flow of 
\[  H_0\rest_M = \frac { \mu v^2 }2 - \frac{\mu^3 }6 \,, \]
with respect to the symplectic form $ \omega\rest_M $. 

 The systems of equations \eqref{eq:t3} and \eqref{eq:t3slow} are 
obtained using the following basic idea: if a Hamilton flow of $ H $, 
with initial condition on a symplectic submanifold, $ M $, 
stays close to $ M $, then the flow is close to the Hamilton 
flow of $ H \rest_M $. 

In our case $ M$ is the manifold of solitons and $ H $ is given 
by \eqref{eq:HV}
\begin{equation}
\label{eq:Heff}
 H_V\rest_M ( a ,v, \gamma, \mu ) =
\frac { \mu v^2 }2 - \frac{\mu^3 }6 +
\frac12 \mu^2 (V  * \sech^2) ( \mu a ) \,, 
\end{equation}
and in particular
\[ H_{ q \delta_0 } \rest_M = 
H_0 \rest_M + 
\frac12 \mu^2 \, \sech^2 ( \mu a ) \,, \ \
H_{ W( h \bullet) } \rest_M = H_0 \rest_M + \frac12 \mu^2 \, W ( h \bullet ) 
* \sech^2 ( \mu \bullet ) \,. \]
The equations \eqref{eq:t3} are simply the equations of the flow of 
$ H_{q \delta_0 } \rest_M $ -- see \S \ref{s:GPH}. The equations of
the flow of $ H_{W ( h \bullet) } \rest_M $ are easily seen to imply
\eqref{eq:t3slow} but some $ h$ corrections are built into the 
classical motion. It would be interesting to see if this provides
improvement of the analysis of \cite{FrSi}. Since our interests
lie in the study of various aspects of the delta impurity we
satisfy ourselves with a numerical experiment which shows that
the improvement is indeed dramatic -- see Fig.\ref{f:drama}.

In either case, 
all of this hinges on the proximity of the flow to $ M$ 
and to show
that we use the Lyapunov function, $ L ( w ) $, introduced in \cite{We} -- see
\S \ref{elf}. Typically, and as is done in \cite{FrSi}, $ L ( w ) $ is bounded
from below so that it controls the norm of $ w $ 
(roughly speaking the expression 
estimated in \eqref{eq:t1} and \eqref{eq:slow}), while
$ (d/dt) L ( w ) $ is estimated from above.
In this paper due to the irregularity of the potential that approach
for upper bounds does not seem to be  applicable but we can estimate $ L (w ) $
directly, controlling the propagation of $ a, v, \gamma $, and $ \mu$
more precisely.

The paper is organized as follows. In \S \ref{rhs} we recall the Hamiltonian
structure of the nonlinear flow of \eqref{eq:nls} and describe the
manifold of solitons. Its identification with the Lie group $ G = 
H_3 \ltimes \RR_+ $, where $ H_3 $ is the Heisenberg group, provides
useful notational shortcuts. In \S \ref{re} we describe
the reparametrized evolution. The starting point there is an 
application of the implicit function theorem and a decomposition of the 
solution into symplectically orthogonal components. That method
has a long tradition in soliton stability and we learned it from \cite{FrSi}. 
In \S \ref{se} we give a self-contained and constructive presentation 
of well known spectral estimates. Weinstein's Lyapunov function is
adapted to our problem in \S \ref{elf}. It is estimated using classical 
energy. The ODE estimates needed for the iteration of our stability 
argument are given in \S \ref{ode} and a stronger version of Theorem \ref{t:1}
is proved in \S \ref{pr}.

Finally, we make comments on the numerics. The computations of 
solutions of \eqref{eq:nls} and \eqref{eq:slow} were done using 
the {\tt FORTRAN} code described in \cite[\S 3]{HMZ2} and written
as part of that project by J. Marzuola. Other computations and 
all the graphics were done using {\tt MATLAB}.

\noindent
{\sc Acknowledgments.} We would like to thank Jeremy Marzuola for allowing 
us the use of his code for NLS computations, and to 
Patrick Kessler and Jon Wilkening for generous help with various
computing issues.
The work of the first author was supported in part by an NSF postdoctoral 
fellowship, and that
of the second second author by an NSF grant DMS-0200732.

\section{The Hamiltonian structure and the manifold of solitons}
\label{rhs}

In this section we recall the well known facts about the
Hamiltonian structure of the nonlinear Schr\"odinger equation.
The manifold of solitons is given as an orbit of a semidirect
product of the Heisenberg group and $ \RR_+ $.

\subsection{Symplectic structure}
\label{s:sy}

Let $V$ be a complex Hilbert space with the inner product 
$ \langle \bullet, \bullet \rangle_V $. For $W$, a totally real 
subspace of $ V $ ($ W \cap i W = \{ 0 \} $), we have 
$V=W+iW \simeq W^2$, and we can consider $W$ and $ V$ 
as real Hilbert spaces. 

As a real Hilbert space $ V $ is equipped with the natural 
inner product or metric
\[ g(X,Y) = \Re \la X,Y \ra_V \,, \] 
and the natural symplectic form
\[ \omega(X,Y) = \Im \la X, Y \ra_V = g(X,iY)\,. \]
In other words $ g$, $ \omega $, and $ J $, multiplication by $ 1/i $
form a compatible triple:
\begin{equation}
\label{eq:nats} 
\omega ( X , Y ) = g ( J X , Y) \,, \ \ g ( X, Y ) = \omega ( X, i Y) \,.
\end{equation}
In terms of $W^2$, we have
$$g(X,Y) = \left< \begin{bmatrix} \Re X \\ \Im X \end{bmatrix}, 
\begin{bmatrix} \Re Y \\ \Im Y \end{bmatrix} \right>_{W^2} = 
\omega \left( 
\begin{bmatrix} \Re X\\ \Im X \end{bmatrix},  J  
\begin{bmatrix} \Re Y\\ \Im Y \end{bmatrix}  \right)$$
and
$$\omega(X,Y) = g \left( J 
\begin{bmatrix} \Re X\\ \Im X \end{bmatrix}, 
\begin{bmatrix} \Re Y\\ \Im Y \end{bmatrix}  \right)$$
where $J$ is the matrix representing  multiplication by $-i$:
$$J = \begin{bmatrix} 0 & I \\ - I & 0 \end{bmatrix}$$
For example, when we consider $V=\mathbb{C}^n$ and $W=\mathbb{R}^n$, then 
$\omega$ is just the standard symplectic form.  

In our work, we  take 
$V=H^1(\mathbb{R}, \mathbb{C})\subset L^2(\mathbb{R},\mathbb{C})$, and on $V$ 
we use the $L^2$ inner product.  The symplectic form $\omega$ is thus
\begin{equation}
\label{eq:omega}
\omega(u,v) = \Im \int u\bar v\,, 
\end{equation}
and the metric $g$ is
$$g(u,v) = \Re \int u\bar v$$
Now we consider Hamiltonians and associated Hamiltonian flows.  Let 
$H:V\to \mathbb{R}$ be a function, our Hamiltonian.  The associated Hamiltonian
vector field is a map $\Xi_H : V\to TV$, 
which means that for a particular point $u\in V$, we have 
$(\Xi_H)_u \in T_uV$. The vector field  $\Xi_H$ is defined by the relation
\begin{equation}
\label{eq:Hamvf} \omega(v , (\Xi_H)_u) = d_uH(v)
\,, \end{equation}
where $v\in T_uV$, and $d_uH:T_uV \to \mathbb{R}$ is defined by
$$d_uH(v) = \frac{d}{ds}\Big|_{s=0} H(u+sv) \,. $$
In the notation of \eqref{eq:nats} if we use $ g $ to define functionals,
$  dH_u ( v ) = g ( v, \nabla H_u ) $, then $ (\Xi_H)_u = J \nabla H_u  $.

If we take $ V = H^1 ( \RR , \CC ) $ with the symplectic form \eqref{eq:omega},
and
$$H(u) = \int \frac14|\partial_x u|^2 - \frac14|u|^4$$
then we can compute
\begin{align*}
d_uH(v) &= \Re \int ( (1/2) \partial_x u \partial_x \bar v - |u|^2u \bar v)\\
&= \Re \int ( - (1/2)\partial_x^2 u - |u|^2u)\bar v \,. 
\end{align*}
Thus, in view of \eqref{eq:nats} and \eqref{eq:Hamvf}, 
$$(\Xi_H)_u = \frac 1 i \left( - \frac12 \partial_x^2u - |u|^2u \right) $$
The flow associated to this vector field (Hamiltonian flow) is
\begin{equation}
\label{eq:Hflow}
\dot u = (\Xi_H)_u =  \frac 1 i  \left( - 
\frac12 \partial_x^2u - |u|^2u \right) 
\,.
\end{equation}

For future reference we state two general lemmas of symplectic 
geometry:
\begin{lem}
\label{l:gen1}
Suppose that $ g : V \rightarrow V $ is a diffeomorphism such 
that $ g^* \omega = \mu ( g ) \omega $, where $ \mu( g ) \in C^\infty ( 
V ; \RR)  $. Then for $ f \in C^\infty ( V , \RR ) $,
\begin{equation} 
\label{eq:lg1}
(g^{-1})_* \Xi_f(g(\rho)) = \frac{1}{\mu(g)}\Xi_{g^*f}(\rho) \,, \ \ 
\rho \in V \,. 
\end{equation}
\end{lem}
\begin{proof}
This is a straightforward generalization of Jacobi's theorem 
which is the case of $ \mu ( g ) \equiv 1 $.
To compute $(g^{-1})_* \Xi_f(g(\rho)) $, we note
\begin{align*}
  \omega_\rho((g^{-1})_*X, \, (g^{-1})_*\Xi_f(g(\rho))) 
&=((g^{-1})^*\omega)_{g(\rho)}(X, \, \Xi_f(g(\rho)))
= \frac{1}{\mu(g)}\omega_{g(\rho)}(X,\Xi_f(g(\rho))) \\
&=\frac{1}{\mu(g)} [ df(g(\rho))] (X) 
= \frac{1}{\mu(g)} [g^*df(\rho)] ((g^{-1})_*X ) \\
&= \frac{1}{\mu(g)}\omega_\rho((g^{-1})_*X, \Xi_{g^*f}(\rho))
\end{align*}
and the lemma follows.
\end{proof}

Suppose that $ f \in C^\infty ( V ; \RR ) $ and that $ df ( \rho_0  ) = 0 $.
Then the Hessian of $ f $ at $ \rho_0 $, 
$ f'' ( \rho_0 ) : T_{\rho_0} V \mapsto T^*_{\rho_0} V $, is well defined. 
The Hamiltonian map $ F : T_{\rho_0} V \rightarrow T_{\rho_0} V $ is given 
by the relation
\begin{equation} 
\label{eq:HamF}
\left[ f''(\rho_0 ) X \right] (Y )  = \omega_{\rho_0} ( Y , F X )
\,. 
\end{equation}

In this notation we have 
\begin{lem}
\label{l:gen2}
Suppose that $ N \subset V $ is a finite dimensional symplectic submanifold of 
$ V $, 
and $ f \in C^\infty ( V , \RR )$ satisfies 
\[  \Xi_f ( \rho ) \in T_\rho N \subset T_\rho V \,, \ \ \rho \in N \,.\]
If at $ \rho_0 \in N $, $ df ( \rho_0) = 0 $, then the Hamiltonian map
defined by \eqref{eq:HamF} satisfies
\[  F ( T_{\rho_0} N ) \subset T_{\rho_0} N  \,. \]
\end{lem}
\begin{proof} Since $ N $ is assumed to be finite dimensional 
we only need to prove the lemma for a finite dimensional $ V$ (any 
particular $ Y \in (T_\rho V)^\perp $ can be a value of a vector field
in a finite dimensional submanifold of $ V $ containing $ N $). 
We can then assume that $ \rho_0 = ( 0 , 0 ) $, and that in local 
coordinates near $ ( 0 , 0 ) $, $ N = \{ ( x, \xi ) \; | \; x''=\xi''=0 \}, $
$ x = ( x' , x'')$ , $ \xi = ( \xi' , \xi'') $, $ \bullet'= 
( \bullet_1 , \cdots, \bullet_k ) $, where $ 2k = \dim N $ (see 
for instance \cite[Theorem 21.2.4]{Hor2}). The conditions of $ f $ mean
that 
\[  d_{x''} f ( x',\xi',0,0) = d_{ \xi''} f ( x',\xi',0,0) = 0 \,, \ \ 
d f ( 0 , 0 ) = 0 \,,\]
where we wrote $ ( x, \xi ) = ( x',\xi',x'',\xi'') $. 
Hence, the Hessian at $ ( 0 , 0 ) $ is given by 
\[  f''( 0 , 0 ) = \begin{bmatrix} f''_{x',\xi'} ( 0 , 0 ) & 0 \\
0 & f''_{x'',\xi''} ( 0 , 0 ) \end{bmatrix} \,. \]
This means that 
\[       \langle f''(\rho_0 ) 
X , Y \rangle = 0 \ \ \forall\; X \in T_\rho N\,, 
\ Y \in ( T_\rho N )^\perp \,.\]
where $ \bullet^\perp $ denotes the symplectic orthogonal. Since 
the Hamiltonian map, $ F $, is defined by $ \langle f''(\rho_0 ) X , Y \rangle 
= \omega ( Y, F Y ) $ this proves the lemma.
\end{proof}

\subsection{Associated symmetries and Noether's theorem}
For completeness we comment on 
the Hamiltonian version of Noether's theorem which states that the
following three statements are equivalent
\begin{gather*}
\Xi_H E \defeq \omega ( \Xi_H , \Xi_E )  =0\,, \\  \text{$E$ is preserved by the Hamiltonian flow of $H$,}
\\ \text{$H$ is preserved by the Hamiltonian flow of $E$.}
\end{gather*}
For example, consider the mass $M = \int |u|^2$.  The associated Hamiltonian 
vector field is $(\Xi_M)_u = iu$.  We compute
$$\omega(\Xi_M,\Xi_H) = 
- \Im \int iu \, \overline{i(\partial_x^2 u + |u|^2u)} =0$$
The flow associated to $\Xi_M$ is $ u \mapsto e^{is} u $, 
which is the phase invariance of 
$H$ and thus solutions to $\partial_t u= i(\partial_x^2u + |u|^2u)$.

Similarly, the time translation, $ u ( x, t ) \mapsto u ( x , t + s ) $
gives the conservation of energy, $ H ( u ) $, 
the space translation, $ u ( x , t ) \mapsto u ( x +y , t  ) $, 
gives the conservation of momentum, $ \Im \int u_x \bar u $. 

\subsection{Manifold of solitons as an orbit of a group}

For $ g = ( a, v, \gamma , \mu ) \in \RR^3 \times \RR_+ $ we define
the following map 
\begin{equation}
\label{eq:repG}    H^1 \ni u \longmapsto g\cdot u \in H^1 \,, \ \ 
(g\cdot u)(x) \defeq e^{i\gamma}e^{iv(x-a)}\mu u(\mu(x-a)) \,. 
\end{equation}
This action gives a group structure on $ \RR^3 \times \RR_+ $
and it is easy to check that this
transformation group is a semidirect product of the Heisenberg group
$ H_3 $ and $ \RR_+ $:
\[ G= H_3\ltimes\mathbb{R}_+ \,, \ \ 
\mu\cdot(a,v,\gamma) = (\frac{a}{\mu} ,\mu v, \gamma) \,.\]
We recall that the Heisenberg group can be identified with the 
group of matrices of the form 
$$\begin{bmatrix}
1 & v & \gamma \\
0 & 1 & a \\
0 & 0 & 1
\end{bmatrix}\,, \ \ a, v , \gamma \in \RR \,,
$$
and that the semidirect product of $ H $ and $ \RR_+ $ is defined
by 
\[ ( h , \mu ) \cdot ( h' , \mu' ) 
= ( h \cdot ( \mu \cdot h' ) , \mu \mu' ) \,, \ \ h , h' \in H \,.\] 
Explicitly, the group law on $ G $ is given by 
\[ (a,v,\gamma,\mu) \cdot (a',v',\gamma',\mu') = (a'',v'',\gamma'',\mu'') \,, 
\] 
where 
\[ \begin{split}
& v'' = v + v'\mu \,, \ \ 
 a'' = a + \frac{a'}{\mu} \,, \ \
 \gamma'' = \gamma + \gamma' + \frac{va'}{\mu} \,, \ \
 \mu'' = \mu \mu'
\end{split} \]

\noindent
{\bf Remark.} As was pointed to us by Bjorn Poonen, the 
group acts faithfully on the 4-dimensional space spanned
by $1$,$v$,$a$,$\gamma$ viewed as functions on the group. This can
be used to see that the group is faithfully represented by 
the group of matrices of the form 
$$\begin{bmatrix}
1 & 0 & 0 & 0 \\ v & \mu & 0 & 0 \\ a & 0 & {1}/{\mu} & 0 \\ \gamma & 0 & 
{v}/{\mu} & 1 
\end{bmatrix} \,, \ \ v, a , \gamma \in \RR \,, \ \mu \in \RR_+ \,,
$$
but we will not use this below.

The action of $ G $ is not symplectic but it is {\em conformally symplectic}
in the sense that 
\begin{equation} 
\label{eq:conf} g^*\omega = \mu(g)\omega
\,, \ \ g = ( h ( g ) , \mu ( g ) ) \,, \ \ \mu ( g ) \in \RR_+ \,,
\end{equation}
as is easily seen from \eqref{eq:omega}.

The Lie algebra of $ G$, denoted by $ {\mathfrak g } $, 
is generated by $ e_1, e_2, e_3 , e_4 $, 
\[ \begin{split} & \exp( t e_1 ) = (t ,0,0,1)  \,, \ \  
\exp( t e_2 ) = (0 ,t,0,1)  \,,
 \\ &  \exp(t e_3 ) = ( 0, 0 , t, 1 ) \,, \ \ \exp ( t e_4 ) 
= ( 0 , 0 , 0 , e^t) \,, \end{split} \]
and the bracket acts as follows:
\begin{equation}
\label{eq:liea}
[e_1,e_4]=e_1, \quad [e_2,e_4] = -e_2, \quad [e_1,e_2]=-e_3, 
\quad [e_3,\bullet ] =0 \,, \end{equation}
so $e_3$ is in the center. 
The infinitesimal representation obtained from \eqref{eq:repG} 
is given by
\begin{equation}
\label{eq:liea1}  e_1 = -\partial_x \,, \ \
e_2 = ix \,,  \ \ e_3 = i \,, \ \  e_4 = \partial_x \cdot x \,. 
\end{equation}
It acts, for instance on $ {\mathcal S}( \RR ) \subset H^1 $, 
and by $ X \in {\mathfrak g} $ we will denote a linear combination of
the operators $ e_j $.

We have the following standard
\begin{lem}
\label{l:stan}
Suppose $ \RR \ni t \mapsto g ( t) $ is a $ C^1 $ function and that
$ u \in {\mathcal S} ( \RR ) $. Then, in the notation of \eqref{eq:repG},
\[  \frac{d}{dt } g(t) \cdot u = g( t ) \cdot ( X ( t ) u ) \,, \]
where $ X( t ) \in {\mathfrak g } $ is given by 
\begin{equation}
\label{eq:lstan}
X ( t ) =   \dot a  ( t) \mu ( t) e_1 +
 \frac {\dot v (t) } { \mu ( t ) } e_2 + 
(\dot \gamma ( t) - \dot a ( t ) v ( t )  ) e_3 + 
\frac  {\dot\mu( t) }{ \mu ( t ) } e_4 \,, 
\end{equation}
where $ g ( t) = ( a ( t ) , v ( t) , \gamma( t ) , \mu ( t) ) $.
\end{lem}
\begin{proof}
We differentiate
\[  g ( t ) \cdot u = \exp ( i \gamma ( t) ) \exp ( - a ( t) \partial_x ) 
\exp ( i v ( t) x ) \exp ( \theta ( \partial_x \cdot x ) ) u \,, \ \
\exp \theta ( t) = \mu ( t ) \,,\]
and note that 
\[ \begin{split}
&\partial_x \exp ( i v  x ) = \exp ( i v x ) ( \partial_x + i v ) \,, \\ 
&\partial_x  \exp ( \theta ( \partial_ x \cdot x ) )   = 
 \exp ( \theta ( \partial_ x \cdot x ) ) e^\theta \partial_x \,, \\
&i x  \exp ( \theta ( \partial_ x \cdot x ) )   = 
 \exp ( \theta ( \partial_ x \cdot x ) ) ( e^{-\theta } i x ) \,, 
\end{split} \]
either by direct computation or using \eqref{eq:liea}. The formula 
\eqref{eq:lstan} follows.
\end{proof}

The manifold of solitons is an orbit of this group, $ G \cdot \eta $, 
to which 
$ \Xi_H $, defined in \eqref{eq:Hamvf}, is tangent. In view of
\eqref{eq:Hflow} that means that 
\[  i \left( \frac12 \partial_x^2 \eta + |\eta|^2 \eta \right) = 
X \cdot \eta \,, \]
for some $ X \in {\mathfrak g } $. The simplest choice is given 
by taking $ X = \lambda i  $, $ 
\lambda \in \RR $, so that $ \eta $ solves
a nonlinear elliptic equation 
$$-\frac12 \eta'' - \eta^3 + {{\lambda}}\eta =0 \,. $$
This has a solution in $ H^1 $ if $ \lambda = \mu^2/2 > 0 $ and it then is
$\eta(x)= \mu \sech(\mu x)$.
We will fix $ \mu = 1$ so that 
$$\eta(x)=  \sech x \,. $$ 
Using Lemma \ref{l:gen1} we can check that $ G\cdot \eta $ is
the {\em only} orbit of $ G $ to which $ \Xi_H $ is tangent. 

We define the submanifold of solitons, $M \subset H_1$,
as the orbit of $\eta$ under $G$, 
$$M = G\cdot \eta \subset H_1$$
and thus we have the identifications
\begin{equation}
\label{eq:idG}
  M = G\cdot \eta \simeq G / \ZZ \,, \ \  T_\eta M = {\mathfrak g} \cdot 
\eta \simeq  {\mathfrak g }  \,. 
\end{equation}
The quotient corresponds to the $\ZZ$-action 
$$ ( a, v, \gamma, \mu)  \mapsto (a, v, \gamma + 2 \pi k , \mu ) \,, \ \ 
 k \in \ZZ $$

\subsection{Symplectic structure on the manifold of solitons}
\label{ssms}
We first compute the symplectic form $ \omega \rest_M $ on 
$ T_\eta M $ using the identification \eqref{eq:idG}:
\[ (\omega\rest_M )_\eta ( e_i , e_j ) = \Im \int (e_i \cdot \eta )(x)
(\overline {e_j \cdot \eta} ) ( x ) \,.\]
Since 
\[ \int \eta^2 ( x ) dx = 2 \,, \ \ \int \eta ( x) \partial_x \eta ( x ) = 0 
\,, \ \  \int \partial_x \eta ( x ) x \eta ( x ) dx = -1 \,, \]
we obtain from \eqref{eq:liea1} that 
\begin{equation}
\label{eq:omij}
  (\omega\rest_M )_\eta ( e_2 , e_1 ) = 1 \,, \ \ (\omega\rest_M)_\eta 
( e_3, e_4 ) = 1 \,, 
\end{equation}
and all the other $   (\omega\rest_M )_\eta ( e_i , e_j ) $'s vanish.
In other words,
\[    (\omega\rest_M )_\eta = (dv \wedge da + d \gamma \wedge d \mu )_
{ ( 0 ,0 , 0 , 1 ) } = ( d ( v da + \gamma d\mu ) )_{ ( 0, 0 , 0 , 1 )} \,. \]

To find an expression for $ \omega\rest_M $ we use \eqref{eq:conf}
and the following elementary
\begin{lem}
\label{l:elem}
If $ \sigma $ is a one form on $ \RR^3 \times \RR_+ $ such that
\[  \sigma_{(0,0,0,1)} =   (  v da + \gamma d\mu  )_{ ( 0, 0 , 0 , 1 )} 
 \,, \ \  g^* \sigma = \mu ( g) \sigma \,, \ \ g \in G \,, 
\] 
then 
\[ \sigma  = \mu v da + \gamma d \mu \,.\]
\end{lem}
\stopthm

We conclude that using the identification \eqref{eq:idG} 
\begin{equation}
\label{E:Momega}
\omega \rest_M  = \mu dv\wedge da + vd\mu\wedge da + d\gamma \wedge d\mu
\end{equation}

Now let $f$ be a function defined on 
$M$, $ f = f ( a, v, \gamma, \mu ) $.
The associated Hamiltonian vectorfield, $\Xi_f$, is defined by 
$$\omega(\cdot, \Xi_f) = 
df = f_ada+f_vdv+f_\mu d\mu + f_\gamma d\gamma \,. $$
Using \eqref{E:Momega} we obtain 
\begin{equation}
\label{eq:XifM}
\Xi_f = \frac{f_v}{\mu} \partial_a 
+ \left( -\frac{f_a}{\mu}-\frac{vf_\gamma}{\mu}\right)\partial_v 
+ f_\gamma \partial_\mu + \left(v\frac{f_v}{\mu}-f_\mu\right)\partial_\gamma
\,. \end{equation} 
The Hamilton flow is obtained by solving
\[ 
\dot v = -\frac{f_a}{\mu}-\frac{vf_\gamma}{\mu}\,, \ \
\dot a = \frac{f_v}{\mu} \,, \ \
\dot \mu = f_\gamma \,, \ \ 
\dot \gamma = v \frac{f_v} \mu - f_\mu \,. \]
The restriction of 
$$ H ( u) = \frac14\int |\partial_x u|^2 - \frac14\int |u|^4$$ 
to $ M $ is given by computing by 
\begin{equation}
\label{eq:ffM}
 f ( a, v, \gamma, \mu ) = H ( g\cdot \eta ) = 
 \frac{\mu v^2}{2}-\frac{\mu^3}{6} \,, \ \ g = ( a, v , \gamma, \mu ) \,. 
\end{equation}
The flow of \eqref{eq:XifM} for this $ f $ describes the evolution of
a soliton.

\subsection{The Gross-Pitaevski Hamiltonian restricted to the 
manifold of solitons}
\label{s:GPH}

We now consider the Gross-Pitaevski Hamiltonian for the delta function 
potential
\begin{equation}
\label{eq:GPH}
H_q ( u ) \defeq \frac14 \int ( |\partial_x u |^2 - |u|^4 ) dx 
+ \frac12 q | u ( 0 ) |^2 \,, 
\end{equation}
and its restriction to $ M = G \cdot \eta $:
\begin{equation}
\label{eq:GPHM} 
{H_q }\rest_M = f ( a , v , \gamma, \mu ) = 
\frac { \mu v^2 }2 - \frac{\mu^3 }6 + 
\frac12 q \mu^2 \sech^2 ( \mu a ) \,. 
\end{equation}
This is obtained from 
\eqref{eq:ffM} and from calculating
\[  \frac12 q | ( g \cdot \eta ) | ( 0 ) = \frac12 q \mu^2 \eta^2 ( - \mu a ) 
= \frac12 q \mu^2 \sech^2 ( \mu a ) \,. \]
The flow of $ (H_q)\rest_M $ can be read off from \eqref{eq:XifM}:
\begin{equation}
\label{eq:GPfl}
 \begin{split}
&\dot v = -\frac{f_a}{\mu} - \frac{v f_\gamma}{\mu} = \mu^2 q \, \sech^2(\mu a) 
\tanh(\mu a)\\
&\dot a = \frac{f_v}{\mu} = v \\
&\dot \mu = f_\gamma = 0\\
&\dot \gamma = v \frac{f_v} \mu - f_\mu = \frac12 v^2 + \frac12 \mu^2 - 
q\mu \, \sech^2(\mu a) 
- \frac12q\mu^2 a \, \sech^2(\mu a) \tanh(\mu a)
\end{split} \end{equation}
This are the same equations as \eqref{eq:t3}. The evolution of $ a $ and 
$ v $ is simply the Hamiltonian evolution of 
$ ( v^2 + q \mu^2 \sech^2 ( \mu a ) )/2 $, $ \mu = \text{const} $. 
The more mysterious evolution of the phase $ \gamma $ is now explained
by \eqref{eq:GPHM}.

\begin{figure}
\begin{center}
\includegraphics[width=6in]{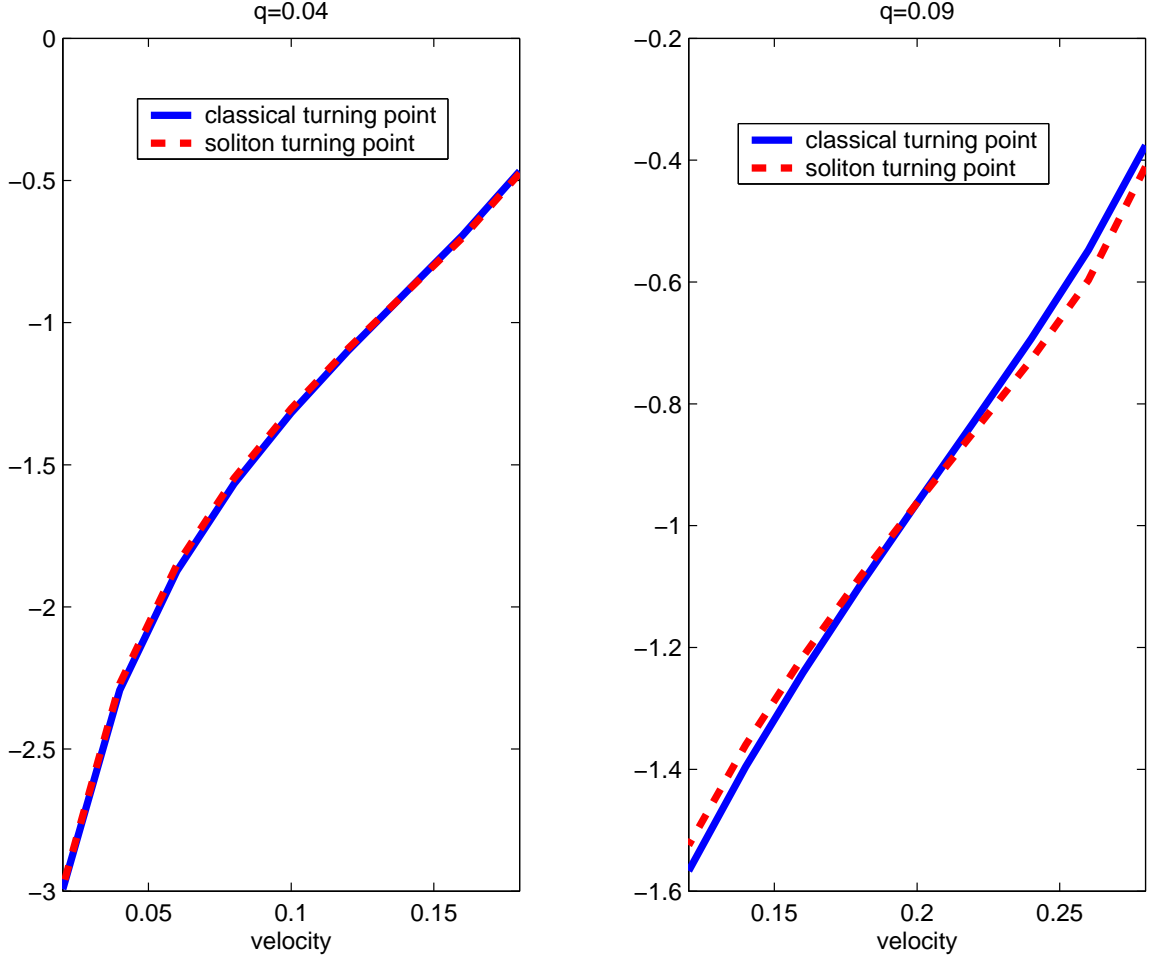}
\end{center}
\caption{
Two plots with $ q = 0.04 $ and $ q=0.09 $, respectively, and $ a_0= -10$.
The blue line is the theoretical prediction of the 
turning point of the soliton, 
$ |a_{\text{turn}}| = 
\sech^{-1} \left( {v}/{\sqrt q}\right)$, and the red dashed line
is the actual soliton turning point. For smaller values of $ q $ the
agreement is outstanding.}
\label{f:turn}
\end{figure}

Since $\mu$ is constant by the third equation, solving this system reduces to 
solving the first two equations. The  
turning position, $a_\text{turn}$, is given by 
$$|a_\text{turn}| = \sech^{-1} \left( \frac{v}{ \sqrt q} \right)$$
and Fig.\ref{f:turn} gives a comparison between  $ a_\text{turn} $ and
the numerically computed turning point of the center of the
soliton.

\section{Reparametrized evolution}
\label{re}

To see the effective dynamics described in \S \ref{s:GPH} we 
write the solution of \eqref{eq:nls} as 
\[  u ( t) = g ( t ) \cdot ( \eta + w ( t )  ) \,,  \ \ 
w ( t) \in H^1 ( \RR, \CC ) \,, \]
where $ w ( t) $ satisfies
\[  \omega ( w ( t) , X \eta ) = 0 \,, \ \ \forall X \in {\mathfrak g}\,.\]
To see that this decomposition is possible, 
initially for small times, we apply the 
following consequence of the implicit function theorem and the nondegeneracy 
of $ \omega\rest_M $ (see \cite[Proposition 5.1]{FrSi} for a more
general statement):
\begin{lem} 
\label{l:FGJS}
For $\Sigma \Subset G/\ZZ$ (where the topology on $G/\ZZ $ is given by 
the identification with  $ \mathbb{R}\times 
\mathbb{R} \times S^1 \times \mathbb{R}_+$) let 
$$U_{\Sigma,\delta} = \{ \, u \in H_1 \, : \, 
\inf_{g\in\Sigma}\|u -g\cdot \eta\|_{H^1}<\delta \} \,. $$
If $ \delta \leq \delta_0 = \delta_0 (\Sigma)$ then for any $ 
u \in U_{\Sigma,\delta}$, there exists a unique 
$  g( u )\in \Sigma $ such that
\[ \omega(g(u)^{-1} \cdot u -  \eta, X \cdot \eta ) = 0 \quad \forall 
 X \in {\mathfrak g } \,. \] 
Moreover, the map $ u  \mapsto g( u )$ is in 
$C^1(U_{\Sigma,\delta},\Sigma)$. 
\end{lem}
\begin{proof}
We define the transformation
\[ F \; : \; H^1 ( \RR , \CC ) \times G \; \longrightarrow \; 
{\mathfrak g}^* \,, \ \ [ F ( u , h ) ]( X)  \defeq 
\omega ( h \cdot u - \eta , X \cdot \eta ) \,.\]
We want to solve 
$ F ( u , h ) = 0 $ for $ h = h ( u ) $ and by the implicit
fuction theorem that follows for $ u $  near $ G \cdot \eta $ if for any 
$ g_0 \in G $ the linear transformation
\[  d_h F ( g_0 \cdot \eta , g_0 ) \; : \; T_{g_0} G \longrightarrow {\mathfrak
g}^* \,, \]
is invertible. Clearly we only need to check it for $ g_0 = e$, that is
that 
$  d_h F ( \eta , e ) \; : \; {\mathfrak g}  \rightarrow {\mathfrak
g}^* \,, $ 
is invertible. But as an element of 
$ {\mathfrak g}^* \otimes {\mathfrak g}^* $, 
$ d_h F( \eta , e ) = (\omega\rest_M)_\eta $, which is nondegenerate.
\end{proof}

For \S\S \ref{s:sy} and \ref{s:GPH} we recall that the equation for 
$ u $ \eqref{eq:nls} can be written as 
\begin{equation} 
\label{eq:utH}
\partial_t u  = \Xi_{H_q} (u ) \,, \ \ 
H_q ( u ) \defeq \frac14 \int ( |\partial_x u |^2 - |u|^4 ) dx 
+ \frac12 q | u ( 0 ) |^2 \,. \end{equation}
Using Lemma \ref{l:FGJS} we define
\begin{equation}
\label{eq:wt}  w(t) = g(t)^{-1} u (t) - \eta \,, \ \ 
g ( t ) \defeq g ( u ( t) ) \,, 
\end{equation}
and we want to to derive an equation for $ w ( t ) $. 

By the chain rule and Lemma \ref{l:stan} 
\[ \begin{split}
\partial_t u ( t) & = \partial_t (  g(t) \cdot ( \eta + w ( t ) ) )
 = g ( t) \cdot ( Y ( t) ( \eta + w ( t ) ) + \partial_t w ( t ) ) \,,
\end{split} \]
\[ Y ( t) \defeq   \dot a  ( t) \mu ( t) e_1  +
 \frac {\dot v (t) } { \mu ( t ) } e_2 + 
(\dot \gamma ( t) - \dot a ( t ) v ( t )  ) e_3 + 
\frac  {\dot\mu( t) }{ \mu ( t ) } e_4 \,, 
 \]
$ g ( t) = ( a ( t ) , v ( t ) , \gamma( t ) , \mu ( t ) ) $.
Combined with \eqref{eq:utH} this gives
\begin{equation}
\label{eq:wtH}
\partial_t w(t) = - Y ( t) \eta - Y ( t ) w + g(t)^{-1}\Xi_{H_q} ( (g(t)
\cdot (\eta 
+ w ( t ) )))  \,.
\end{equation}
To make this more explicit
we apply Lemma \ref{l:gen1} to see that 
$$ g ( t)^{-1} \Xi_{H_q} g( t )  = \frac{1}{\mu(t)}\Xi_{g( t) ^*H_q}$$
(since the action of $ g ( t ) $ is linear on $ H^1 $, $ g(t)^{-1} $
and $ ( g ( t)^{-1})_* $ are identified). We compute
\begin{equation}
\label{eq:gstH}
\begin{split}
g^*H_q ( \tilde u ) &= \frac14 \int ( \mu |\partial_x(e^{ixv}\tilde u ( \mu x ) )|^2 
- \mu^4 | \tilde u (\mu x ) |^4 ) dx  + \frac12 q \mu^2 
|\tilde u(-\mu a)|^2 \\
&= \frac14 \int ( \mu^3  |\partial_x \tilde u ( x ) |^2  - 2 v \mu^2 
\Im \partial_x \tilde u ( x ) \overline{ {\tilde u} ( x )}
+ \Re  {v^2} \mu |\tilde u( x ) |^2- 
\mu^3 | \tilde u  ( x ) |^4 ) dx \\ 
& \ \ \ \ \ \ \ \ \ \ 
+ \frac12 \mu^2 q |\tilde u(-\mu a)|^2 \,, 
\end{split}
\end{equation}
so that
$$ \frac 1 { \mu ( g) }\Xi_{g^*H_q}(\tilde u) = \frac{1}{i} \left( 
 -\frac{\mu^2}{2}\tilde u_{xx} +  v \mu 
\tilde u_x- \mu^2 | \tilde u |^2 \tilde u +\frac{v^2}{2} \tilde u 
+ \mu q\delta( \bullet +\mu a)\tilde u \right) \,. $$
For us $ \tilde u ( t) = \eta + w ( t) $ and we expand the nonlinear 
term 
$$|\eta+w|^2(\eta+w) = \eta^3 + \underbrace{2\eta^2w + \eta^2\bar
w}_{\text{linear}} + \underbrace{2|w|^2\eta + \eta w^2}_{\text{quadratic}} +
\underbrace{|w|^2w}_{\text{cubic}} $$

Inserting this in \eqref{eq:wtH} gives
\begin{lem}
\label{l:wt}
If $ w ( t) $ is given by \eqref{eq:wt} then 
\begin{equation}
\label{eq:lw1}
\partial_t w = X ( t) w + X ( t) \eta -
i\mu^2\mathcal{L}w +i\mu^2\mathcal{N}w 
-iq \mu \delta_0({x} + {\mu}a)\eta -iq \mu \delta_0({x} + {\mu}a)w\,,
\end{equation}
where $ X ( t) \in {\mathfrak g }$ is given by 
\begin{equation}
\label{eq:lw2}
X ( t) \defeq  
\left(-\mu \dot a + v \mu \right)e_1
- \frac{\dot v}{\mu} e_2 + 
\left(-\dot \gamma + v \dot a - \frac{v^2}2 +\frac{\mu^2}{2} \right)e_3 
  -\frac{\dot \mu}{\mu} e_4  \,, \ \ 
\end{equation}
and
$$\mathcal{L}w = - \frac12 \partial_x^2 w - 2\eta^2 w - \eta^2 
\bar w + \frac12 w  \,, \ \ \mathcal{N}w =  2|w|^2\eta + \eta w^2  + |w|^2w\,.$$
\end{lem}
\stopthm

We now want to estimate the coefficients of $ X (t ) $ in \eqref{eq:lw1}
using the symplectic orthogonality of $ Y \eta $, $ Y \in {\mathfrak g} $
and $ w $.  For that we define
\[  P \; : \; {\mathcal S}'( \RR, \CC ) \longrightarrow {\mathfrak g} \]
as the unique linear map satisfying
\[  \omega ( u - P(u)\eta, Y \eta ) = 0 \ \ \forall \, Y \in {\mathfrak g} 
\,. \]
We will need the following

\begin{lem}
\label{l:P}
Let $ \| \bullet \| $ be a norm on $ \mathfrak g $ obtained by 
using the standard $ \RR^4 $ norm in the basis given by \eqref{eq:liea1}.
Then for $ w \in H^1 $, and $ Y \in {\mathfrak g} $,
\[ \begin{split} & \| P ( Y w ) \| \leq C \| Y\| \|w \|_{L^2} \,, \\ 
& \| P ( i \mathcal{N} u ) \| \leq C \|w\|_{L^2}^2 \left( 1 + 
\| w \|_{H^1}^{\frac12} \| w \|_{L^2}^{\frac12} \right) \,, \\
& \| P ( ( i\delta_0({x} - x_0 )w ) \| \leq C 
\| w \|_{H^1}^{\frac12} \| w \|_{L^2}^{\frac12} \,, 
\end{split} \]
with the constant independent of $ x_0 $.
\end{lem}
\begin{proof}
We start with an explicit expression for $ P $ which follows from 
\eqref{eq:omij}:
\begin{equation}
\label{eq:Pexp}
\begin{split}
 P & = \sum_{j=1}^{4} e_j P_j \,, \ \ P_j \; : \; {\mathcal S}' 
\longrightarrow \RR 
\,, \\
P_1 ( u ) & = - \omega ( u , e_2 \eta ) =  \Re \int u ( x ) x \eta ( x ) dx \,, \\
P_2 ( u ) & =  \omega ( u , e_1 \eta ) = - \Im \int u ( x) \partial_x\eta( x ) dx \,, \\
P_3 ( u ) & = \omega ( u , e_4 \eta ) = \Im \int u ( x ) \partial_x ( x \eta ( x)) dx \,, \\
P_4 ( u ) & = - \omega ( u , e_3 \eta) = \Re \int u ( x ) \eta ( x ) dx \,.
\end{split}
\end{equation}
We now recall that $ \| u \|^2_{L^\infty ( \RR ) } \leq C \| u \|_{ L^2 ( \RR)}
\| u \|_{H^1 ( \RR ) } $ and the estimates follow.
\end{proof}

Since $ P w = P w_t = 0 $, \eqref{eq:lw1} gives
\begin{equation}
\label{eq:X1} \begin{split}
X_1 ( t ) & \defeq X ( t) - q \mu P ( i \delta_0 ( \bullet + a \mu) \eta ) 
\\ 
& =
- P ( X ( t) w) + \mu^2 P (i \mathcal{L}w ) - \mu^2 P ( i \mathcal{N}w ) 
+ q \mu P ( i\delta_0({x} + {\mu}a)w ) \,.
\end{split}
\end{equation}
Since $ \mathcal{L}$ is the Hessian of $ \mathcal{E} $, given in 
\eqref{eq:defE}, at the 
critical point $ \eta $, and $ \Xi_{\mathcal E} $ is tangent to $ M$, 
Lemma \ref{l:gen2} (or a direct computation) shows that
\[  P ( i \mathcal{L} w ) = 0 \,, \]
and hence that term can be dropped from the right hand side.
We can then use 
Lemma \ref{l:P} to obtain
\begin{prop}
\label{p:1}
Suppose that $ w ( t ) $ is given in Lemma \ref{l:wt} and
that $ X_1 ( t ) $ is given by \eqref{eq:X1}. Then 
\begin{equation}
\label{eq:Xt} 
\| X_1 ( t)\| \leq C q \| w \|_{H^1} + C ( \|w\|_{L^2}^2 + \| w \|_{H^1}^3) 
\,. 
\end{equation}
\end{prop}
\stopthm

Finally we interpret the coefficients of $ X_1 ( t ) $. First 
we use \eqref{eq:Pexp} to see that
\[ P ( i \delta_0 ( \bullet + a \mu) \eta ) = 
\frac12 \partial_{x} ( \eta^2 ) (  a \mu ) e_2 + \left( \eta^2 ( a \mu ) + 
 \frac12 a \mu \partial_x (\eta^2 ) ( a \mu )\right) e_3 \,.\]
Then we combine this with 
\eqref{eq:lw2} and  \eqref{eq:X1} to obtain
\begin{equation}
\label{eq:X12t} \begin{split} X_1 ( t) & =   
\left(-\mu \dot a + v \mu \right)e_1 + 
\left( - \frac12 q \mu 
\partial_{x} ( \eta^2 ) (  a \mu ) - \frac{\dot v}{\mu} \right)
e_2 \\ 
& \ \ \ + 
\left( - \mu q \eta^2 ( a \mu ) 
- \frac12 q a \mu^2 \partial_x (\eta^2 ) ( a \mu )
-\dot \gamma + v \dot a - \frac{v^2}2 +\frac{\mu^2}{2} \right)e_3 
  -\frac{\dot \mu}{\mu} e_4  \,. \end{split} 
\end{equation}
We now see that
\[  X_1 ( t ) = 0  \ \Longleftrightarrow \ \text{equations \eqref{eq:GPfl}
hold.} \]

\section{Spectral estimates}
\label{se}

In this section we will recall the now standard estimates on the
operator $ {\mathcal L} $ which arises as Hessian of $ {\mathcal E} $ 
at $ \eta $:
\[ {\mathcal L} w = - \frac12 \partial_x^2 w - 2 \eta^2 w - \eta^2 \bar w 
+ \frac12 w \,, \]
or 
\[ {\mathcal L} w = \begin{bmatrix} L_+ &  0 \\  0 & L_- \end{bmatrix} 
\begin{bmatrix} \Re w \\ \Im w \end{bmatrix} \,, 
 \ \ L_\pm = - \frac 12 \partial_x^2 - ( 2 \pm 1 ) \eta^2 + \frac12 \,.
\]
In our special case we can be more precise than in the general case
(see \cite{We}, and also \cite[Appendix D]{FrSi}). The 
self-adjoint operators
$ L_\pm $ belong the class of Schr\"odinger operators with 
{\em P\"oschl-Teller} potentials
and their spectra can be explicitly computed using hypergeometric 
functions -- see for instance \cite[Appendix]{GuZw}. This gives
\[ \sigma( L_- ) = \{ 0 \} \cup [1/2 , \infty ) \,, \ \
\sigma(L_+) = \{ 0, - 3/2\} \cup [1/2, \infty ) \,. \] 
The eigenfuctions can computed by the same method but a straightforward
verification is sufficient to see that 
\[ L_- \eta = 0 \,, \ \ L_+ (\partial_x \eta) = 0 \,, \ \ 
L_+ (\eta^2) = - \frac32 \eta^2 \,.\]

We now have 
\begin{prop}
\label{p:coer}
Suppose that for every $ X \in {\mathfrak g}$
\[ \omega ( w , X \cdot \eta ) = 0 \,, \ \ w \in H^1 ( \RR, \CC )  \,. \]
Then, with $ \langle w , v \rangle \defeq \Re \int w \overline v $ on $ 
H^1 ( \RR , \CC ) $ (considered as a real Hilbert space),
\begin{equation}
\label{eq:coer}
\langle {\mathcal L} w , w \rangle \geq \rho_0 \| w\|^2_{L^2} \,, \ \ 
\rho_0 =  \frac{9}{2(12+\pi^2)} \simeq 0.2058 \,.
\end{equation}
\end{prop}
We need 
the following elementary
\begin{lem}
\label{l:lina}
Let $ V $ be a real vector space with an inner product 
$ \langle \bullet, \bullet \rangle $, and let $ L $ be 
a symmetric operator on $ V$.
Suppose that for $ v_0 , v_1 \in V $, $ \| v_j \| = 1 $,  we have 
\begin{gather}
\label{eq:lina1} 
\begin{gathered}
 L v_0 = - c_0  v_0 \,, \ \ c_0 \geq 0 \,, \ \ 
\langle v_0 , v_1 \rangle^2 = c_2  \,, \\ 
\langle w , v_0 \rangle = 0 \ \Longrightarrow \ \langle L w , w \rangle 
\geq c_1 \| w \|^2  \,, \ \ c_1 \geq 0 \,.
\end{gathered}
\end{gather}
Then 
\begin{equation}
\label{eq:lina2}
\langle v , v_1 \rangle =0 \ \Longrightarrow \ 
\langle L v , v \rangle \geq c_3 \| v \|^2 \,, \ \ 
c_3 \defeq  c_1 c_2- {c_0} ( 1 - c_2 )   \,. 
\end{equation}
\end{lem}
\begin{proof}
For reader's convenience we present 
the straightforward argument in which we can assume that $ 0 < c_2 < 1 $.
For $ v \in V $ we write $ v = \alpha v_0 + w $, $ \langle v_0, w \rangle 
= 0 $. The condition $ \langle v , v_1 \rangle = 0 $ gives 
\begin{equation}
\label{eq:alpha}   \alpha^2 =  \frac{1}{c_2} \langle w, v_1 \rangle^2 = 
 \frac{1}{c_2} \langle w , v_1 - c_2^{\frac12} v_0 \rangle^2 \leq 
\frac{1 - c_2 }{c_2} \| w \|^2 \,. 
\end{equation}
Hence 
\[
\begin{split} 
\langle L v , v \rangle & \geq c_1 \|w \|^2 - c_0 \alpha^2 \\
& \geq c_1 \delta  \|w\|^2 + \left( c_1 ( 1 - \delta) 
\frac{c_2} { 1 - c_2 } - c_0 \right) \alpha^2 \\
& = (  c_1 c_2- {c_0} ( 1 - c_2 )  ) \| v \|^2 \,,
\end{split}
\]
if we choose $ c_1 \delta =  (  c_1 c_2 - {c_0} ( 1 - c_2 )  ) $. 
\end{proof}

\noindent
{\em Proof of Proposition \ref{p:coer}:} The assumption means
that 
\[ \Im \int w \left\{ \begin{array}{l} \ \ i \eta \\
\ \partial_x \eta \\ \ i x \eta \\ \partial_x (x \eta) \end{array}
\right\} dx = 0 \,. \]
Working with real and imaginary 
parts the proof reduces to lower bounds on $ L_\pm $:
\begin{gather}
\label{eq:Lpm}
\begin{gathered}
\langle v , \eta \rangle = \langle v , x \eta \rangle = 0 
\ \Longrightarrow \ \langle L_+ v , v \rangle \geq \rho_0 \|v \|^2_{L^2}
\,, \\
\langle v , \partial_x 
\eta \rangle = \langle v , ( x\partial_x +1 )\eta \rangle = 0 
\ \Longrightarrow \ \langle L_- v , v \rangle \geq \rho_0 \|v \|^2_{L^2}
\,,
\end{gathered}
\end{gather}
where now $ v \in H^1 ( \RR; \RR ) $. 
Noting that 
\[ \langle \eta , \partial_x \eta \rangle = \langle x \eta , \eta^2 \rangle
= \langle \partial_x \eta, \eta^2 \rangle = 0 \]
we can apply Lemma \ref{l:lina}
in the following three cases:
\begin{gather*}
 V = (\partial_x \eta)^\perp \cap H^2 ( \RR, \RR ) \,, \ \ v_0 = \frac{\sqrt{3}}2 \eta^2 \,, 
\ \ v_1 = \frac{1}{ \sqrt 2} \eta \,, \ \ L = L_+ \\
c^1_0 = \frac32 \,, \ \ c^1_1 = \frac12 \,, \ \ c^1_2 = \frac{3 \pi^2}{32}
\,, \\
V = ( \eta^2)^\perp  \cap H^2 ( \RR, \RR ) \,, \ \ v_0 = \frac{\sqrt{3}}{\sqrt 2} 
\partial_x \eta \,, \ \ 
v_1 = \frac{\sqrt{6}}{\pi} x \eta \,, \ \ L = L_+ \\
c^2_0 = 0 \,, \ \ c^2_1 = \frac12 \,, \ \ c_2^2 = \frac{9 }{\pi^2} \,,
\\
V = H^2 ( \RR , \RR ) \,, \ \ 
v_0 = \frac{1}{ \sqrt 2} \eta \,, \ \ v_1 = \frac{2\sqrt 2}{\sqrt{12+\pi^2}} \partial_x ( x \eta) \,,  \ \ L = L_- \,, \\
c^3_0=0 \,, \ \ c_1^3 = \frac12\,, \ \ c_2^3 = \frac{9}{12+\pi^2}  \,.
\end{gather*}
Here we used
\begin{gather*}
\int_\RR \sech^2 ( x ) dx = 2 \,, \ \ \int_\RR \sech^4 ( x) 
 dx = \frac43  \,,  \ \ 
\int_\RR \sech^3 ( x) dx = \frac{\pi}{2} \,, \ \ 
\int_\RR x^2 \sech^2 (x ) dx = \frac{\pi^2}{6}
\,, \\ 
\int_\RR \tanh^2 ( x ) \sech^2 ( x ) dx = \frac23 \,, \ \ 
\int_\RR ( \partial_x ( x \sech( x ) ))^2 dx = \frac{1}{18}(12+\pi^2) \,.
\end{gather*}
It follows that we can take
\[ \rho_0 = \min_{j=1,2,3}( c_1^jc_2^j - c_0^j ( 1 - c_2^j) )
= \min \left( \frac{3 \pi^2 }{16} - \frac32 , \frac {9}{2\pi^2} , 
\frac{9 }{2 (12 + \pi^2)} \right) =  \frac{9 }{2 (12 + \pi^2)} \,,\]
completing the proof.
\stopthm

Proposition \ref{p:coer} gives a slightly stronger statement:
\[ 
\begin{split}
\langle {\mathcal L} w , w \rangle & \geq 
( 1 - \delta ) \langle {\mathcal L} w , w \rangle + 
\delta \rho_0 \| w\|^2_{L^2}  \\ 
& \geq ( 1 - \delta ) \left( \frac12 \| \partial_x w \|^2 - 
\frac{5}{2} \| w \|^2 \right) + \delta \rho_0 \|w \|^2_{L^2} \\
& \geq \frac{2 \rho_0} { 5 + 2 \rho_0 } \| \partial_x w \|^2  
\simeq 0.0760 \| \partial_x w \|^2 \,, \ \ 
\delta = \frac{5}{5 + 2 \rho_0 } \,.
\end{split} \] 
In addition, 
\begin{equation}
\label{eq:lowerL}
\langle {\mathcal L} w , w \rangle  \geq 
 \frac{2 \rho_0  } { 7 + 2 \rho_0 }  \| w \|^2_{H^1} 
\simeq 0.0555 \| w \|_{H^1}^2 \,.
\end{equation}

\medskip
\noindent
{\bf Remark.} The smallness of these constants gives a possible
explanation of the size of $ q$'s for which the asymptotic
result agrees with numerical simulations. The implicit constants
in \S \ref{elf} are closely related to the constants above.

\renewcommand\thefootnote{\ddag}%

\section{Estimates on the Lyapunov function}
\label{elf}

Suppose $u=u(x,t)$ solves \eqref{eq:nls} 
with\footnote{The symbol $\ll 1$ means smaller than an \textit{absolute} 
positive constant, i.e.\ one independent of all parameters in this problem.} 
$|q|\ll 1$ and initial data 
\begin{equation}
\label{E:init_cond}
u_0(x) = e^{ixv_0}\eta(x-a_0), \qquad |v_0| \ll 1
\end{equation}
Let $T>0$ be the maximal time such that on $[0,T]$, the smallness condition 
$\delta \leq \delta_0$ in Lemma \ref{l:FGJS} is met.  From Lemma \ref{l:FGJS}, 
obtain the $C^1$ parameters $\mu=\mu(t)$, $\gamma=\gamma(t)$, $v=v(t)$, 
$a=a(t)$ satisfying the symplectic orthogonality conditions stated there.  
Let $\tilde u=\tilde u(x,t)$ be defined by
\begin{equation}
\label{E:psi}
u(x,t)= g(t) \cdot \tilde u ( x , t ) \defeq 
 e^{i\gamma}e^{ixv}\mu \tilde u(\mu (x-a),t) \, ,
\end{equation}
and let
$$w(x,t)=\tilde u(x,t)-\eta(x) \, .$$ 

The Lyapunov function of \cite{We} and \cite{FrSi} is given by 
\begin{equation}
\label{eq:Lya} 
L ( w ) \defeq {\mathcal E} ( \eta + w ) - {\mathcal E} ( \eta ) \,.
\end{equation}
The lower bound on $ L ( w ) $ follows from the spectral estimates of 
\S \ref{se}, and in particular from \eqref{eq:lowerL}. 
For the upper bound we will use the conservation of $ H_q ( u ) $
and its relation to $ {\mathcal E}( \eta + w ) $.

For future reference we state the following crucial consequence of  
the orthogonality conditions on $ w $, and in particular of the 
condition that $  \Im \int  i \eta \bar w  =  \Re \int \bar w \eta =0 $:
\begin{lem}
\label{l:gath}
Suppose that for every $ X \in {\mathfrak g}$
\[ \omega ( w , X \cdot \eta ) = 0 \,, \ \ w \in H^1 ( \RR, \CC )  \,. \]
Then 
\begin{equation}
\label{eq:lga1}
\|w\|_{L^2}^2 = \frac2 \mu  (1 - \mu )\,, \ \ 
\end{equation}
\end{lem}
\begin{proof}
We first compute
\[ \| \eta + w \|_{L^2}^2 = \| g^{-1} u  \|_{L^2}^2  = 
\frac{1}{ \mu(g)} \| u  \|_{L^2}^2 = 
\frac{2}{ \mu(g)} \,,\]
where we used the conservation of the $ L^2 $ norm. As noted before the 
statement of the lemma $ \Re \la w , \eta \ra = 0 $ and hence
\[ \| \eta + w \|_{L^2}^2  = 2 + \| w \|_{L^2}^2 \,, \]
from which the conclusion follows. 
\end{proof}

As a consequence, we can 
dispense with $\mu$ in the estimates 
of Proposition \ref{p:1},
and we reformulate it as 
\begin{prop}
\label{C:nomus}
Suppose $1-\mu \ll 1$ and $|q|\leq 1$.  Then
\begin{gather*}
|v-\dot a| + |\dot v + q \partial_x \eta^2(a)/2| +
|-q \eta^2(a) - q a \partial_x \eta^2(a)/2  -\dot \gamma + v^2/2 + 1/2 |  \\
\leq C(|q| \|w\|_{H^1}^2 + \|w\|_{H^1}^2 + \|w\|_{H^1}^3) \,.
\end{gather*}
\end{prop}
\begin{proof} We use \eqref{eq:lga1} in \eqref{eq:Xt}.
For example, 
\begin{align*}
\Big|\frac12 q \partial_x \eta^2(a) + \dot v \Big| &\leq \mu \Big| \frac12 
\frac{q}{\mu} \partial_x \eta^2 (a) + \frac{\dot v}\mu \Big| \\
&\leq \mu \Big| \frac12 q \mu \partial_x \eta^2(a\mu) + \frac{\dot v}{\mu} 
\Big| + c|q||1-\mu|\\
&\leq 2 \Big| \frac12 q \mu \partial_x \eta^2(a\mu) + \frac{\dot v}{\mu} \Big| 
+ c|q| \|w\|_{L^2}^2
\end{align*}
We also use the estimate for $|v-\dot a|$ to replace $v\dot a$ by $v^2$ in the equation for $\dot \gamma$.
\end{proof}

We adopt the following notational convention:  
denote the initial (time $t=0$) 
configuration of the system by $0$-subscripts -- $u_0=u(0)$, $w_0=w(0)$,  and 
$a_0=a(0)$, $v_0=v(0)$, $\mu_0=\mu(0)$, $\gamma_0=\gamma(0)$.  Similarly, 
denote the configuration of the system at some fixed time $t_i$ by 
$i$-subscripts.  Finally, the configuration at any arbitrary time $t$ we 
denote without subscripts -- $w=w(t)$, $u=u(t)$ and $a=a(t)$, $v=v(t)$, 
$\mu=\mu(t)$, $\gamma=\gamma(t)$.

With this notation we now state
\begin{lem}
\label{L:elf}
Suppose $\mu_0=1$ and $w_0=0$ (equivalently, suppose \eqref{E:init_cond} 
holds), and suppose that $T>0$ is the maximal time for which the smallness 
condition in Lemma \ref{l:FGJS} holds.  Suppose that for an interval of time 
$[t_i,t_{i+1}]\subset [0,T]$, the following conditions hold 
\begin{equation}
\label{E:hypoth}
\begin{gathered}
0\leq 1-\mu \ll 1, \quad \max_{t_i\leq s \leq t_{i+1}} |v(s)|\ll 1, \quad 
\|w_i\|_{H^1} \leq 1 \,, \\
|q||t_{i+1}-t_i| \ll 1, \quad |t_{i+1}-t_i| \max_{t_i\leq s \leq t_{i+1}} 
|v(s)| \ll 1\,. 
\end{gathered}
\end{equation}
Then there is an absolute constant $c_*>1$ such that
$$\sup_{t_i\leq s\leq t_{i+1}} \|w(s)\|_{H^1}^2 \leq c_*\|w_i\|_{H^1}^2 
+ c_*|q|^2 \,. $$
\end{lem}
We remark that the inequality, $ 0 \leq 1 - \mu $, in \eqref{E:hypoth} 
is not an assumption but follows
from Lemma \ref{l:gath}.

The main result of this section is the following consequence of this:
\begin{prop}
\label{C:elf_iterated}
Suppose $\mu_0=1$ and $w_0=0$, and suppose that $T>0$ is the maximal time for 
which the smallness condition in Lemma \ref{l:FGJS} holds.  Let 
$$n \leq \frac{\delta \log(1/|q|)}{\log c_*} -1$$
and suppose there is a partition of the time axis
$$0=t_0<t_1<\cdots < t_n\leq T$$
such that on each subinterval $[t_i,t_{i+1}]$, \eqref{E:hypoth} in Lemma 
\ref{L:elf} holds. Then,
$$\sup_{0\leq s \leq t_n} \|w(s)\|_{H^1}^2 \leq |q|^{2-\delta}$$  
\end{prop}

\begin{proof}[Proof of Lemma \ref{L:elf}]
We start by noting that in the argument that follows, 
we will not use any information about 
$w$ or the parameters $\mu$, $\gamma$, $a$, and $v$ for times $0<t< t_i$; only
that $\mu_0=1$ and $w_0=0$. 

We will conveniently reexpress $ L ( w ) $ given by \eqref{eq:Lya} using
the conserved Hamiltonian and mass. Since $ u = g \cdot \tilde u $, 
we recall \eqref{eq:gstH} 
to obtain:
\begin{equation}
\label{E:Ham_u}
H_q(u)= g^* H_q ( \tilde u ) = 
\begin{aligned}[t]
&\frac14 \mu v^2 \int |\tilde u|^2 + \frac12 \mu^2 v \Im \int 
\partial_x \tilde u \, \bar{\tilde u} 
+ \frac14 \mu^3 \int |\partial_x \tilde u|^2 \\
&-\frac14\mu^3 \int |\tilde u|^4 + \frac12 q\mu^2 |\tilde u(-\mu a,t)|^2
\end{aligned}
\end{equation}
The expression for the mass,
$M(u) = \int |u|^2$, becomes
$ M(u) = \mu \int |\tilde u|^2 $.
Using this and \eqref{E:Ham_u}, we obtain
$$
\mathcal{E}(\tilde u) =
\begin{aligned}[t]
&\frac{1}{\mu^3}H_q(u) + \frac{1}{4\mu}M(u) - \frac{v^2}{4\mu^3}M(u) 
-\frac{v}{2\mu}\Im \int \bar{\tilde u} \partial_x \tilde u 
- \frac{q}{2\mu}|\tilde u(-\mu a)|^2
\end{aligned}
$$
Now substitute $\tilde u=\eta+w$ and use the orthogonality condition 
$\Im \int w \partial_x \eta =0$ to obtain
\begin{equation}
\label{E:Lyap-final}
\mathcal{E}(\eta+w) =
\begin{aligned}[t]
&\frac{1}{\mu^3}H_q(u) + \frac{1}{4\mu}M(u) - \Big( \frac{v^2}{4\mu^3}M(u) 
+ \frac{q}{2\mu}\eta(-\mu a)^2\Big)\\
&-\frac{v}{2\mu}\Im \int \bar w \partial_x w 
- \frac{q}{\mu}\eta(-\mu a)\Re w(-\mu a) - \frac{q}{2\mu}|w(-\mu a)|^2
\end{aligned}
\end{equation}
Note that the classical energy term (with the $\mu$ terms dropped) 
$$ E( u ) \defeq \frac 14 v^2M(u) + \frac12q\eta( a)^2\,,  $$
has appeared in this expression.  
Evaluate \eqref{E:Lyap-final} at $t=t_i$ to obtain
\begin{equation}
\label{E:Lyap-init}
\mathcal{E}(\eta+w_i) =
\begin{aligned}[t]
&\frac{1}{\mu_i^3}H_q(u) + \frac{1}{4\mu_i}M(u) 
- \Big( \frac{v_i^2}{4\mu_i^3}M(u) + \frac{q}{2\mu_i}\eta(-\mu_i a_i)^2\Big)\\
&-\frac{v_i}{2\mu_i}\Im \int \bar w_i \partial_x w_i 
- \frac{q}{\mu_i}\eta(-\mu_i a_i)\Re w(-\mu_i a_i) 
- \frac{q}{2\mu_i}|w(-\mu_i a_i)|^2
\end{aligned}
\end{equation}
By taking the difference of the right hand sides of 
\eqref{E:Lyap-final} and \eqref{E:Lyap-init}, 
we obtain
\begin{equation}
\label{E:Lyap}
\begin{aligned}
\indentalign \mathcal{E}(\eta + w) - \mathcal{E}(\eta) \\
&=
\begin{aligned}[t]
&\Big( \frac{1}{\mu^3} - \frac{1}{\mu_i^3} \Big)H_q(u) 
+ \frac14\Big(\frac1\mu-\frac{1}{\mu_i}\Big)M(u) \\
&-\Big( \frac{v^2}{4\mu^3}M(u) + \frac{q}{2\mu}\eta(-\mu a)^2\Big)
+\Big( \frac{v_i^2}{4\mu_i^3}M(u) + \frac{q}{2\mu_i}\eta(-\mu_i a_i)^2\Big)\\
&-\frac{v}{2\mu}\Im\int\bar w\partial_x w 
- \frac{q}{\mu}\eta(-\mu a)\Re w(-\mu a) - \frac{q}{2\mu}|w(-\mu a)|^2 \\
&+ \frac{v_i}{2\mu_i}\Im \int \bar w_i\partial_x \bar w_i 
+ \frac{q}{\mu_i}\eta(-\mu_ia_i)\Re w_i(-\mu_ia_i) 
+ \frac{q}{2\mu_i}|w_i(-\mu_ia_i)|^2\\
&+(\mathcal{E}(\eta+w_i)-\mathcal{E}(\eta))
\end{aligned}\\
&=\text{I}+\text{II}+\text{III}+\text{IV}+\text{V}
\end{aligned}
\end{equation}
where each line has been labeled by a Roman numeral.  From the spectral 
estimate Proposition \ref{p:coer} (see \eqref{eq:lowerL}), we have
\begin{equation}
\label{E:spec_est}
c_1\|w\|_{H^1}^2 - \|w\|_{H^1}^3 -\frac14\|w\|_{H^1}^4 \leq \mathcal{E}(\eta+w)-\mathcal{E}(\eta)
\end{equation}
We next estimate the right-hand side of \eqref{E:Lyap}, line by line.  
For $t_i \leq t\leq t_{i+1}$, let
$$\epsilon(t)^2 = \sup_{t_i\leq s \leq t} \|w(s)\|_{H^1}^2$$

\noindent \textit{Estimate of the 1st line of \eqref{E:Lyap}}.  By the 
assumption $w_0=0$ and $\mu_0=1$, we have
\begin{equation}
\label{E:mass}
M(u)=M(\eta)=2
\end{equation}
and
\begin{equation}
\label{E:energy}
H_q(u) = -\frac16 + \frac12v_0^2+\frac{q}{2}\eta^2(a_0)
\end{equation}
By substituting \eqref{E:mass} and \eqref{E:energy} into Term I, we obtain 
$\text{I} = \text{I}_a+\text{I}_b$, where
$$\text{I}_a = -\frac16\Big( \frac{1}{\mu^3}-\frac{1}{\mu_i^3}\Big)
+\frac12\Big(\frac{1}{\mu}-\frac{1}{\mu_i}\Big)$$
and
$$\text{I}_b = \Big(\frac{1}{\mu^3}-\frac{1}{\mu_i^3}\Big)\Big(\frac12v_0^2
+\frac{|q|}2\eta^2(a_0)\Big)$$
Inserting \eqref{eq:lga1} in 
Term $\text{I}_a$, gives
\begin{align*}
\text{I}_a&=\frac16\Big(\frac1\mu-\frac1{\mu_i}\Big)\Big( 3- \frac1{\mu^2}
-\frac1{\mu\mu_i}-\frac1{\mu_i^2}\Big)\\
&=\frac16\Big(\frac1\mu-\frac1{\mu_i}\Big)\Big[\Big( 1- \frac1{\mu^2}\Big)
+ \Big(1-\frac1{\mu\mu_i}\Big)+\Big(1-\frac1{\mu_i^2}\Big)\Big]\\
&= 
\begin{aligned}[t]
\frac16\Big(-\frac12\|w\|_{L^2}^2 & +\frac12\|w_i\|_{L^2}^2\Big)
\Big[  \Big(\frac12\big(1+\frac1\mu\big)\|w\|_{L^2}^2\Big) \\
&+ \Big(\frac12\|w\|_{L^2}^2+\frac1{2\mu}\|w_i\|_{L^2}^2\Big)
+ \Big(\frac12\big(1+\frac1\mu_i\big)\|w_i\|_{L^2}^2\Big)\Big]
\end{aligned}
\end{align*}
and thus
$$|\text{I}_a| \leq \frac16(\|w\|_{L^2}^2+\|w_i\|_{L^2}^2)^2$$
For Term $\text{I}_b$, we have
$$ \text{I}_b = \Big( \frac1\mu - \frac{1}{\mu_i}\Big)\Big( \frac1{\mu^2}
+\frac1{\mu\mu_i}+\frac1{\mu_i^2}\Big)\Big(\frac12v_0^2
+\frac{|q|}2\eta^2(a_0)\Big)$$
and thus
$$|\text{I}_b| \leq \frac34(v_0^2+|q|) (\|w\|^2+\|w_i\|^2)$$
Collecting these estimates, we obtain
\begin{equation}
\label{E:line1}
|\text{I}| \leq \epsilon^4 + 2(v_0^2+|q|) \epsilon^2
\end{equation}

\medskip
\noindent
{\bf Remark:} This direct calculation is in fact the consequence
of $ d {\mathcal E}_\eta = 0 $. We are using 
\[  {\mathcal E} ( \mu\cdot \eta ) - {\mathcal E} ( \eta ) = 
{\mathcal O} ( ( 1 -  \mu )^2 ) \,, \]
which follows from 
\[ \partial_\mu {\mathcal E}( \mu \cdot \eta ) \rest_{ \mu = 1 } = 0 \,. \]

\medskip

\noindent \textit{Estimate of the 2nd line of \eqref{E:Lyap} (classical 
energies)}.
We compute
\begin{align*}
\partial_t \Big( \frac{v^2}{2} + \frac{q}{2}\eta^2(a) \Big) &= v\dot v 
+ \frac12 q\partial_x \eta^2(a)\dot a\\
&=  \Big( \dot v + \frac12 q\partial_x\eta^2(a)\Big)v 
+ \frac12q\partial_x\eta^2(a)(\dot a-v)
\end{align*}
and thus by Proposition \ref{C:nomus},
$$ \Big| \partial_t \Big( \frac{v^2}{2} + \frac{q}{2}\eta^2(a) \Big) 
\Big| \leq c(\|w\|_{H^1}^2 + |q|\|w\|_{H^1} + \|w\|_{H^1}^3)(|v|+|q|)$$
By the fundamental theorem of calculus,
$$
 \Big|\Big( \frac{v^2}{2} + \frac{q}{2}\eta^2(a) \Big)- \Big( \frac{v_i^2}{2} 
+ \frac{q}{2}\eta^2(a_i) \Big)\Big| 
\leq c(\epsilon^2+|q|\epsilon +\epsilon^3)\big( |t-t_i|\max_{t_i\leq s\leq t} 
|v(s)| + |q||t-t_i|\big) 
$$
As in the proof of Proposition \ref{C:nomus}, we can install $\mu$'s in this 
expression using \eqref{eq:lga1} to obtain
\begin{equation}
\label{E:line2}
 | \text{II}| \leq c(\epsilon^2+|q|\epsilon)
\big( |t-t_i|\max_{t_i\leq s\leq t} |v(s)| + |q||t-t_i| + v^2 + |q|\big) 
\end{equation}

\noindent \textit{Estimate of the 3rd and 4th lines of \eqref{E:Lyap}}.
By the Cauchy-Schwarz inequality and the Sobolev embedding theorem,
$$|\text{III}| \leq |v| \|w\|_{H^1}^2 + |q| \|w\|_{H^1} 
+ |q| \|w\|_{H^1}^2$$
Similarly, 
$$|\text{IV}| \leq |v| \|w_i\|_{H^1}^2 + |q| \|w_i\|_{H^1} 
+ |q| \|w_i\|_{H^1}^2$$
and thus
\begin{equation}
\label{E:line34}
|\text{III}|+|\text{IV}| \leq 2(|v|+|q|)\epsilon^2 + 2|q|\epsilon
\end{equation}

\noindent \textit{Estimate of the 5th line of \eqref{E:Lyap}}.  By definition 
of $\mathcal{E}$, we have  
\begin{equation}
\label{E:Lpert}
\mathcal{E}(\eta+w_i) = \frac14\int |\partial_x \eta + \partial_x w_i|^2 
- \frac14\int |\eta + w_i|^4 + \frac14\int|\eta+w_i|^2
\end{equation}
Substitute into \eqref{E:Lpert} the three expansions:
\begin{align*}
&|\partial_x \eta+\partial_xw_i|^2 = |\partial_x\eta|^2 
+ 2\Re \partial_x \eta \, \partial_x w_i + |\partial_x w_i|^2\\
&|\eta+w_i|^4 = \eta^4 + 4\Re \eta^3 w_i + 2\eta^2(2(\Re w_i)^2+|w_i|^2) 
+ 4\eta(\Re w_i)|w_i|^2 + |w_i|^4\\
&|\eta+w_i|^2 = \eta^2 + 2\eta\Re w_i + |w_i|^2
\end{align*}
and observe that the linear terms cancel since $\eta$ solves $-\frac12\eta 
+ \frac12\eta'' + \eta^3=0$.  Thus, we obtain the estimate
\begin{equation}
\label{E:line5}
|\text{V}| \leq 8\|w_i\|_{H^1}^2+4\|w_i\|_{H^1}^3+\|w_i\|_{H^1}^4 
\leq 10\|w_i\|_{H^1}^2
\end{equation}

This completes the line-by-line estimation of the right-hand side of 
\eqref{E:Lyap}.  By combining \eqref{E:spec_est}, and the estimates 
\eqref{E:line1},\eqref{E:line2},\eqref{E:line34},\eqref{E:line5} for the 
right-hand side of \eqref{E:Lyap}, we obtain
$$c_1\epsilon^2 \leq 
\begin{aligned}[t]
&\epsilon^3 + \frac14\epsilon^4 
+ c(\epsilon^2+|q|\epsilon+\epsilon^3)(|t-t_i|\max_{t_i\leq s\leq t}|v(s)| 
  + |q||t-t_i| + v^2 + |q|) \\
& + [\epsilon^4 + 2(v_0^2+|q|)\epsilon^2] + [2(|v|+|q|)\epsilon^2+2|q|\epsilon] + 10\|w_i\|_{H_1}^2
\end{aligned}
$$
By hypothesis, every $\epsilon^2$ term on the right side has a small 
coefficient, and thus can be absorbed on the left side.  Therefore, we obtain
$$\epsilon^2 \leq c(|q|\epsilon+\|w_i\|_{H^1}^2)$$
By applying the Peter-Paul inequality $|q|\epsilon\leq \frac12c|q|^2 
+ \frac{\epsilon^2}{2c}$, we obtain the desired estimate.
\end{proof}

\begin{proof}[Proof of Proposition \ref{C:elf_iterated}]
Now let
$$\epsilon^2(t) = \sup_{0\leq s \leq t} \|w(s)\|_{H^1}^2$$
On the first interval $[0,t_1]$, we apply Lemma \ref{L:elf} with $i=0$, and 
since $w_0=0$, we obtain
$$\epsilon(t_1)^2 \leq c_* |q|^2$$
On the second interval $[t_1,t_2]$, we apply Lemma \ref{L:elf} with $i=1$, and 
since $\|w_1\|_{H^1}^2 \leq c_* |q|^2$, we obtain
$$\epsilon(t_2)^2 \leq (c_* + c_*^2)|q|^2$$
We continue, and after the $n$ applications, we obtain
$$\epsilon(t_n)^2 \leq c_* \left( \sum_{j=0}^{n-1}c_*^j \right) |q|^2 
= c_* \left( \frac{c_*^n-1}{c_*-1} \right) |q|^2 \leq c_*^{n+1}|q|^2$$
Since we want $c_*^{n+1}q^2 \leq |q|^{2-\delta}$, we require
$$n+1 \leq \frac{\delta \log(1/|q|)}{\log c_*}$$
\end{proof}

\section{ODE analysis}
\label{ode}

The assumptions of Lemma \ref{L:elf} involve estimates on 
$ v ( s ) $. To control these we use Proposition \ref{C:nomus} 
and ODE estimates which we present in this section.

\begin{lem}
\label{L:ODEcompare}
Suppose $q$ is a constant, $|q|\ll 1$,  and $a=a(t)$, $v=v(t)$, 
$\epsilon_1=\epsilon_1(t)$, $\epsilon_2=\epsilon_2(t)$ are $C^1$ real-valued 
functions.  Suppose $f:\mathbb{R}\to\mathbb{R}$ is a $C^2$ mapping such that 
$|f|$ and $|f'|$ are uniformly bounded.  Suppose that on $[0,T]$, 
\begin{equation}
\label{E:ODE}
\left\{
\begin{aligned}
&\dot a = v + \epsilon_1\\
&\dot v = qf(a) + \epsilon_2
\end{aligned}
\right., \qquad
\begin{aligned}
&a(0)=a_0\\
&v(0)=v_0
\end{aligned}
\end{equation}
Let $\bar a=\bar a(t)$ and $\bar v=\bar v(t)$ be the $C^1$ real-valued 
functions satisfying the exact equations
$$
\left\{
\begin{aligned}
&\dot{\bar a} = \bar v \\
&\dot{\bar v} = qf(\bar a) 
\end{aligned}
\right., \qquad
\begin{aligned}
&\bar a(0)=a_0\\
&\bar v(0)=v_0
\end{aligned}
$$
with the same initial data.  Suppose that on $[0,T]$, we have 
$|\epsilon_j| \leq |q|^{2-\delta}$ for $j=1,2$.  Then provided 
$T\leq \delta |q|^{-1/2}\log (1/|q|)$, we have on $[0,T]$ the estimates
$$|a-\bar a|\leq |q|^{1-2\delta}\log(1/ |q|), 
\qquad |v-\bar v| \leq |q|^{\frac32-2\delta}\log(1/|q|)$$
\end{lem}

Before proceeding to the proof, we recall some basic tools.

\noindent\textit{Gronwall estimate}.  Suppose $b=b(t)$ and $w=w(t)$ are $C^1$ 
real-valued functions, $q$ is a constant, and $(b,w)$ satisfy the differential 
inequality:
\begin{equation}
\label{E:diff_ineq}
\left\{
\begin{aligned}
&|\dot b| \leq |w| \\
&|\dot w| \leq |q| |b|
\end{aligned}
\right. ,
\qquad
\begin{aligned}
&b(0)=b_0\\
&w(0)=w_0
\end{aligned}
\end{equation}
Let $x(t)= |q|^{1/2}b(|q|^{-1/2}t)$, $y(t)=w(|q|^{-1/2}t)$.  Then 
$$
\left\{
\begin{aligned}
&|\dot x| \leq |y| \\
&|\dot y| \leq |x|
\end{aligned}
\right. ,
\qquad
\begin{aligned}
&x(0)=x_0=|q|^{1/2}b_0\\
&y(0)=y_0=w_0
\end{aligned}
$$
Let $z(t)=x^2+y^2$.  Then $|\dot z| = |2x\dot x + 2y\dot y| \leq 2|x||y| 
+ 2|x||y| \leq 2(x^2+y^2) = 2z$, and hence $z(t) \leq z(0)e^{2t}$.
Thus
$$
\begin{aligned}
&|x(t)| \leq \sqrt 2\max(|x_0|,|y_0|) \exp(t)\\
&|y(t)| \leq \sqrt 2\max(|x_0|,|y_0|) \exp(t)
\end{aligned}
$$
Converting from $(x,y)$ back to $(b,w)$, we obtain the Gronwall estimate
\begin{equation}
\label{E:Gron}
\begin{aligned}
&|b(t)| \leq \sqrt 2\max(|q|^{1/2}|b_0|,|w_0|)
\frac{\exp(|q|^{1/2}t)}{|q|^{1/2}}\\
&|w(t)| \leq \sqrt 2\max(|q|^{1/2}|b_0|,|w_0|)\exp(|q|^{1/2}t)
\end{aligned}
\end{equation}

\noindent\textit{Duhamel's formula}.
For a two-vector function $X(t): \mathbb{R} \to \mathbb{R}^2$, a two-vector 
$X_0\in \mathbb{R}^2$, and a $2\times 2$ matrix function 
$A(t):\mathbb{R}\to (2\times 2\text{ matrices})$, let $X(t)=S(t,t')X_0$ denote 
the solution to the ODE system $\dot X(t) = A(t)X(t)$ with $X(t')=X_0$.  In 
other words, $\frac{d}{dt} S(t,t')X_0 = A(t)S(t,t')X_0$ and $S(t',t')X_0=X_0$. 
Then, for a given two-vector function $F(t):\mathbb{R}\to \mathbb{R}^2$, the 
solution to the inhomogeneous ODE system  
\begin{equation}
\label{E:inhomODE}
\dot X(t) = A(t)X(t) + F(t)
\end{equation}
with initial condition $X(0)=0$ is given by Duhamel's formula
\begin{equation}
\label{E:Duhamel}
X(t) = \int_0^t S(t,t')F(t')dt'
\end{equation}

\begin{proof}[Proof of Lemma \ref{L:ODEcompare}]
Let $\tilde a= a-\bar a$ and $\tilde v = v-\bar v$; these perturbative 
functions satisfy
$$
\left\{
\begin{aligned}
&\dot{\tilde a} = \tilde v + \epsilon_1\\
&\dot{\tilde v} =  q g \tilde a + \epsilon_2
\end{aligned}
\right., \qquad
\begin{aligned}
&\tilde a(0)=0\\
&\tilde v(0)=0
\end{aligned}
$$
where $g=g(t)$ is given by
$$g=\left\{
\begin{aligned}
&\frac{f( a)-f(\bar a)}{a-\bar a} & \text{if }\bar a \neq a\\
& f'(a) &\text{if }a=\bar a
\end{aligned}
\right.
$$
which is $C^1$ (in particular, uniformly bounded).  Set 
$$A(t) = \begin{bmatrix} 0 & 1 \\ qg(t) & 0 \end{bmatrix}, 
\quad F(t) = \begin{bmatrix} \epsilon_1(t) \\ \epsilon_2(t) \end{bmatrix}, 
\quad X(t)=\begin{bmatrix} \tilde a(t) \\ \tilde v(t) \end{bmatrix}$$
in \eqref{E:inhomODE}, and appeal to Duhamel's formula \eqref{E:Duhamel} to 
obtain
\begin{equation}
\label{E:Duhamel2}
\begin{bmatrix}
\tilde a(t) \\ \tilde v(t) 
\end{bmatrix}
= \int_0^t S(t,t') \begin{bmatrix} \epsilon_1(t') \\ \epsilon_2(t') 
\end{bmatrix} \, dt'
\end{equation}
Apply the Gronwall estimate \eqref{E:Gron} with
$$\begin{bmatrix} b(t) \\ w(t) \end{bmatrix} = S(t+t',t')\begin{bmatrix} 
\epsilon_1(t') \\ \epsilon_2(t') \end{bmatrix},  \quad \begin{bmatrix} b_0 
\\ w_0 \end{bmatrix} = \begin{bmatrix} \epsilon_1(t') \\ 
\epsilon_2(t') \end{bmatrix} $$
to conclude that 
$$\left| S(t,t') \begin{bmatrix} \epsilon_1(t') \\ 
\epsilon_2(t') \end{bmatrix} \right| \leq \sqrt 2 \begin{bmatrix} 
|q|^{-1/2}\exp(|q|^{1/2}(t-t')) \\ \exp(|q|^{1/2}(t-t')) \end{bmatrix} 
\max(|q|^{1/2}|\epsilon_1(t')|,|\epsilon_2(t')|)$$
Feed this into \eqref{E:Duhamel2} to obtain that on $[0,T]$
\begin{align*}
&|\tilde a(t)| \leq \sqrt 2 \, T\frac{\exp(|q|^{1/2}T)}{|q|^{1/2}} 
\sup_{0\leq s\leq T}\max(|q|^{1/2}|\epsilon_1(s)|,|\epsilon_2(s)|)\\
&|\tilde v(t)| \leq \sqrt 2 \, T\exp(|q|^{1/2}T) 
\sup_{0\leq s\leq T}\max(|q|^{1/2}|\epsilon_1(s)|,|\epsilon_2(s)|)
\end{align*}
Taking $T\leq \delta|q|^{-1/2}\log(1/|q|)$, we obtain the claimed bounds.
\end{proof}

\section{Main theorem and proof}
\label{pr}
Here we put all the components together and give a stronger
version of Theorem \ref{t:1}. The basic procedure is the
iteration of Lemmas \ref{L:elf} and \ref{L:ODEcompare} which
can roughly be described as follows: if the conditions \eqref{E:hypoth}
hold, and the initial condition satisfies $ \| w_i \|_{H^1} \leq
|q|^{1-\delta}$, say, then on the interval $ [t_i, t_{i+1} ]$,
$ \| w ( t ) \|_{ H^1}  \leq  2| q|^{1-\delta } $. That means
that the evolution of the parameters $ g ( t) \in G $ is close
to the evolution using the effective Hamiltonian, in the way that
makes Lemma \ref{L:ODEcompare} applicable. But that gives us
a lower bound on $ t_{i+1}$ for which \eqref{E:hypoth} holds on 
$ [ t_i , t_{i+1} ] $, 
closing the bootstrap loop.

More precisely, we have
\begin{thm}
\label{t:2}
Suppose $|q|\ll 1$ and $|v_0| \ll 1$.  Let $u$ solve 
$$i\partial_t u + \partial_x^2 u - q\delta_0(x)u + |u|^2 u =0$$
with initial data $u_0(x) $ 
satisfying 
\[  \| u_0 -  e^{i \bullet v_0} \eta( \bullet -a_0) \|_{H^1} 
\leq C | q | \,.\]    
Then, for times 
$0\leq t \leq \delta (v_0^2+|q|)^{-1/2}\log(1/|q|)$, the smallness condition in Lemma 
\ref{l:FGJS} is met, and thus there are $C^1$ parameters $\mu$, $v$, $\gamma$,
 $a$ satisfying the symplectic orthogonality conditions stated there.  
Furthermore, we have
$$\| u - \mu e^{ixv}e^{i\gamma} \eta( \mu (x-a)) \|_{H^1} 
\leq c|q|^{1-\frac12 \delta}$$
Moreover, if $\bar a$, $\bar v$, $\bar\gamma$ solve the ODE system
\begin{equation}
\label{E:approxODE}
\dot{\bar a} = \bar v \,, \ \
\dot{\bar v} = -\frac12 q\partial_x \eta^2(\bar a) \,, \ \
\dot{\bar \gamma} = \frac12\bar v^2 + \frac12 -q\eta^2(\bar a) 
+ \frac12q\bar a \partial_x \eta^2(\bar a) \,.
\end{equation}
with initial data $(a_0,v_0,0)$, then 
$$
 |a-\bar a| \leq c|q|^{1-3\delta} \,, \ \  |\gamma - \bar \gamma| 
+  |v-\bar v| \leq c|q|^{\frac32-3\delta} \,, \ \ 
 |\mu -1 | \leq c|q|^{2-\delta}\,. 
$$
\end{thm}

\begin{proof}
The equations \eqref{E:approxODE} imply the conservation of energy
$$
\frac12\bar v^2 + \frac12 q\eta^2(\bar a) = \frac12v_0^2+\frac12q \eta^2(a_0)
$$
from which we obtain the bound
\begin{equation}
\label{E:vbarbd}
|\bar v| \leq \sqrt{v_0^2+2|q|} \,. 
\end{equation}
Let 
$$\epsilon(t)^2 = \sup_{0\leq s \leq t} \|w(s)\|_{H^1}^2 \,. $$ 
By Proposition \ref{C:nomus},
\begin{equation}
\label{E:est1}
|\dot a- v| + |\dot v + \frac12q \partial_x\eta^2(a)| \leq 
c_0 ( 
 |q|\|w\|_{H^1}+ \|w\|_{H^1}^2  + \|w\|_{H^1}^3 ) \,. 
\end{equation}

Let $t_1$ with $T\geq t_1>0$ be the maximal time for which 
the assumptions of Lemma \ref{L:elf} \eqref{E:hypoth}
hold with $i=0$.  
Then by Proposition \ref{C:elf_iterated} with $n=1$, we have 
$\epsilon^2(t_1) \leq |q|^{2-\delta}$. The estimate \eqref{E:est1} 
implies \eqref{E:ODE} in Lemma \ref{L:ODEcompare} 
for $ t \in [0,t_1]$, with $f(a)=-\partial_x\eta(a)/2$.   By Lemma 
\ref{L:ODEcompare} and \eqref{E:vbarbd}, we have 
$$\max_{0\leq s\leq t_1} |v(s)| \leq 2 \sqrt{v_0^2+2|q|}\,. $$  
Reviewing \eqref{E:hypoth}, we now see that 
$$T\geq t_1 \geq c_4(v_0^2+2|q|)^{-1/2}\,, $$
where $c_4$ depends only on the implicit absolute constant in \eqref{E:hypoth}.  

Now let $t_2$ with $T\geq t_2>t_1$ be the maximum time such that 
\eqref{E:hypoth} holds with $i=1$.  Then by Proposition \ref{C:elf_iterated} 
with $n=2$, we have 
$$\epsilon^2(t_2) \leq |q|^{2-\delta}\,.$$
  By \eqref{E:est1}, 
we have that \eqref{E:ODE} in Lemma  \ref{L:ODEcompare} holds on $[0,t_2]$.   
By Lemma \ref{L:ODEcompare} and \eqref{E:vbarbd}, we have 
$$ \max_{0\leq s\leq t_2} |v(s)| \leq  2 \sqrt{v_0^2+2|q|} \,. $$
Reviewing \eqref{E:hypoth}, we now see that 
$$ | t_2-t_1| \geq c_4(v_0^2+2|q|)^{-1/2} \,, $$
with the same $c_4$ as in the previous paragraph.

Continue until the $n$th step is reached, where 
$$ n  =  \frac{\delta \log(1/|q|)}{\log c_*} -1 \,, $$ 
which is the most allowed in 
Proposition \ref{C:elf_iterated}.  But now we know that 
$$ 
T\geq t_n \geq c\delta (v_0^2+2|q|)^{-1/2}\log(1/|q|)\,, $$
 and that on $[0,t_n]$, 
$$|a-\bar a| \leq |q|^{1-2\delta}\log(1/|q|), \quad |v-\bar v| 
\leq |q|^{\frac32-2\delta}\log(1/|q|)$$
We also have from Proposition \ref{C:nomus},
$$\Big|\dot \gamma - \Big(\frac12v^2 + \frac12 - q \eta^2( a) 
+ \frac12q a\partial_x\eta(a)\Big) \Big| \leq 2\|w\|_{H^1}^2 + |q|\|w\|_{H^1}$$
Subtracting the equations for $\dot \gamma$ and $\dot{\bar \gamma}$ and using 
that $\|w\| \leq |q|^{2-\delta}$, we obtain
\begin{align*}
|\dot \gamma - \dot{\bar \gamma}| &\leq |v^2-\bar v^2| 
+ |q| |\eta^2(a)-\eta^2(\bar a)| + |q| |a-\bar a| \eta^2(\bar a) 
+ |q||\bar a| |\partial_x^2 \eta(a) - \partial_x^2 \eta (\bar a)| \\
&\leq \left( |q|^{1/2}|q|^{\frac32-2\delta} + |q| |q|^{1-2\delta}
+ |q||q|^{1-2\delta}\right)\log(1/|q|) + |q|^{3-4\delta}\log^2(1/|q|) \\
&\leq |q|^{2-2\delta}\log(1/|q|)
\end{align*}
Since we restrict to times $t \leq \delta|q|^{-1/2}\log(1/|q|)$, we integrate to
obtain
$|\gamma - \bar \gamma| \leq |q|^{\frac32-3\delta}$.
\end{proof}

\medskip
\noindent
{\bf Remark.} There remains the case of initial velocities, 
$ v_0$, which are not small. When $ |q| \rightarrow 0 $
and $ v_0 > 0 $ is fixed, the dynamics is not interesting and the 
solution can be approximated by the solution with $ q =0 $, that
is by the propagating soliton \eqref{eq:frees}. The proof of 
that follows from the arguments of \cite[\S 3.1]{HMZ1}.

\end{document}